\documentclass[12pt,a4paper]{amsart}
\usepackage[vmargin=80pt,hmargin=72pt]{geometry}
\usepackage[T1]{fontenc}
\usepackage[utf8]{inputenc}
\usepackage[english]{babel}
\usepackage{amsfonts,amsmath,amsthm}
\usepackage[bookmarks, colorlinks, breaklinks]{hyperref}
\usepackage[all]{xy}
\usepackage{graphicx}
\usepackage{suffix}
\newtheorem{theorem}{Theorem}[section]

\newtheorem{lemma}[theorem]{Lemma}
\newtheorem{proposition}[theorem]{Proposition}
\theoremstyle{definition}
\newtheorem{definition}[theorem]{Definition}

\theoremstyle{remark}
\newtheorem{example}[theorem]{Example}

\numberwithin{equation}{section}
\newcommand{\rbra}{\left(}       \newcommand{\rket}{\right)}
\newcommand{\sbra}{\left[}       \newcommand{\sket}{\right]}
\newcommand{\cbra}{\left\{}      \newcommand{\cket}{\right\}}
\newcommand{\abra}{\left\langle} \newcommand{\aket}{\right\rangle}
\newcommand{\nbra}{\left.}       \newcommand{\nket}{\right.}
\newcommand{\HH}{\mathcal{H}}
\newcommand{\LL}{\mathcal{L}}
\newcommand{\MM}{\mathcal{M}}
\newcommand{\RR}{\mathbb{R}}
\newcommand{\To}{\longrightarrow}
\newcommand{\Mapsto}{\longmapsto}
\newcommand{\tensor}{\otimes}

\newcommand{\horizontal}{\mathsf{h}}
\newcommand{\Vertical}{\mathrm{Vert}}
\newcommand{\vertical}{\mathsf{v}}

\newcommand{\hor}{\horizontal}
\renewcommand{\Vert}{\Vertical} %use \| instead
\renewcommand{\vert}{\vertical} %use  | instead

\newcommand{\Christoffel}{\Gamma}
\newcommand{\chr}{\Christoffel}

\newcommand{\Smooth}{\mathcal{C}^\infty}
\newcommand{\Sections}{\mathrm{Sec}}
\newcommand{\Fields}{\mathfrak{X}}
\newcommand{\Forms}{\Omega}
\newcommand{\Linear}{\mathrm{Lin}}
\newcommand{\Lin}{\Linear}
\newcommand{\Affine}{\mathrm{Aff}}
\newcommand{\Aff}{\Affine}
\newcommand{\id}{\mathop{\mathrm{id}}\nolimits}
\newcommand{\ex}{\mathop{\mathit{ex}}\nolimits}
\newcommand{\pr}{\mathop{\mathit{pr}}\nolimits}
\newcommand{\leg}{\mathop{\mathrm{leg}}\nolimits}
\newcommand{\Leg}{\mathop{\mathrm{Leg}}\nolimits}
\newcommand{\pairing}{\abra\,\cdot\,,\cdot\,\aket}
\newcommand{\diff}{d}

\newcommand{\du}{\diff u}
\newcommand{\dx}{\diff x}
\newcommand{\dy}{\diff y}
\newcommand{\dmx}{\diff^mx}

\newcommand{\dmxi}{\diff^{m-1}x_i}
\newcommand{\dmxj}{\diff^{m-1}x_j}

\newcommand{\dd}[2][{}]{\frac{\diff#1}{\diff#2}} %\mathchoice?
\WithSuffix\newcommand\dd*[2][{}]{\diff#1/\diff#2} %\ensuremath?

\newcommand{\vardd}[2][{}]{\frac{\delta#1}{\delta#2}} %\mathchoice?
\WithSuffix\newcommand\vardd*[2][{}]{\delta#1/\delta#2} %\ensuremath?
\newcommand{\pp}[2][{}]{\frac{\partial#1}{\partial#2}} %\mathchoice?
\WithSuffix\newcommand\pp*[2][{}]{\partial#1/\partial#2} %\ensuremath?
\newcommand{\ie}{\emph{i.e.}}
\newcommand{\refer}{\emph{ref.}}
\newcommand{\quand}{\quad\textrm{and}\quad}

\newcommand{\secref}[1]{\ref{#1}}%\S

\title[]{Classical field theories of first order and Lagrangian submanifolds of premultisymplectic manifolds}

\author[C. M. Campos]{Cédric M. Campos}
\address{Dept. Matemática Fundamental\\
         Universidad de La Laguna, ULL\\
         Avda. Astrofísico Fco. Sánchez\\
         38206 La Laguna, Tenerife (Spain)}
\email{cedricmc@ull.es}
\author[E. Guzmán]{Elisa Guzmán}
\email{eguzman@ull.es}
\author[J.C. Marrero]{Juan Carlos Marrero}
\email{jcmarrer@ull.es}

\hypersetup{
  linkcolor=blue,
  citecolor=blue,
  filecolor=blue,
  urlcolor=MidnightBlue,
  pdftitle={Tulczjew's triple for Classical Field Theory},
  pdfauthor={C. M. Campos, E. Guzmán, J. C. Marrero},
  pdfsubject={MSC2010 70S05 (70H03, 70H05, 53D12)},
  pdfkeywords={field theory, multisymplectic structure, Lagrangian submanifold, Tulczyjew's triple, Euler-Lagrange equation, Hamilton-De Donder-Weyl equation},
  pdfcreator={GNU Emacs 23.3.1},
}

\begin{document}
\begin{abstract}
A description of classical field theories of first order in terms of Lagrangian submanifolds of premultisymplectic manifolds is presented. For this purpose, a Tulczyjew's triple associated with a fibration  is discussed. The triple is adapted to the extended Hamiltonian formalism. Using this triple, we prove that Euler-Lagrange and Hamilton-De Donder-Weyl equations are the local equations defining Lagrangian submanifolds of a premultisymplectic manifold.
\end{abstract}

\maketitle
\tableofcontents

% INTRODUCTION -----------------------------------------------------------------
\section{Introduction}
It is well-known that the Lagrangian formulation of time-independent Mechanics may be developed using the geometry of the tangent bundle $TQ$ of the configuration space $Q$ and the Hamiltonian formulation may be given using the canonical symplectic structure $\Omega$ of the cotangent bundle $T^*Q$ (see, for instance, \cite{AbrhMrsd78}). But Lagrangian and Hamiltonian Mechanics may also be formulated in terms of Lagrangian submanifolds of symplectic manifolds \cite{Tlcz76a,Tlcz76b} (see also \cite{LnRdrg89}). In fact, the Euler-Lagrange equations for a Lagrangian function on $TQ$ and the Hamilton equations for a Hamiltonian function on $T^*Q$ are just the local equations defining Lagrangian submanifolds of $T(T^*Q)$. Here, the symplectic structure on $T(T^*Q)$ is just the complete lift $\Omega^c$ of $\Omega$. Moreover, in this construction (the so-called Tulczyjew's triple for Classical Mechanics),  an important role is played by the canonical involution of the double tangent bundle $T(TQ)$ and by the isomorphism between $T(T^*Q)$ and $T^*(T^*Q)$ induced by the canonical symplectic structure of $T^*Q$.

For time-dependent Mechanics, the situation is more complicated. Anyway, in the restricted (respectively, extended) formalism of Lagrangian and Hamiltonian time-\-de\-pen\-dent Mechanics may be formulated in terms of Lagrangian submanifolds of Poisson (respectively, presymplectic) manifolds. In this case, we have a fibration $\pi: E\to\RR$, where the total space $E$ is the configuration manifold, which for instance may be $E=\RR\times Q$. As in the time-independent situation, an important role is played by the canonical involution of $T(TE)$ and the canonical isomorphism between $T(T^*E)$ and $T^*(T^*E)$. In fact, using these tools, we prove in our recent paper \cite{GzmMrr10} that, in the restricted formalism, the Euler-Lagrange and Hamilton equations are the local equations defining Lagrangian submanifolds of the Poisson manifold $J^1(\pi_1^*)$.  Here, $\pi^*_1 \colon \Vert^*(\pi) \to \RR$ is the canonical fibration over $\RR$ of the dual bundle of the vertical bundle $\Vert(\pi)$ of $\pi$. We also prove that, in the extended formalism, the dynamical equations are just the local equations defining Lagrangian submanifolds of the presymplectic manifold $J^1\tilde\pi_E$, where $\tilde\pi_E\colon T^*E\to \RR$ is the projection of $T^*E$ onto the real line. Other different descriptions of time-dependent Mechanics in the Lagrangian submanifold setting have been proposed by several authors (see \cite{GrbwUrbn07, IglMrrPdrSs06, LnLcmb89, LnMrr93}).

% Classical Mechanics is a special type of a classical field theory. In classical field theories of first order the Lagrangian function is a real function $L$ on the 1-jet bundle $J^1\pi$, where $ \pi\colon E\to M$ is a fibration (Classical Mechanics corresponds to the particular case when $M$ is the real line). A solution of the Euler-Lagrange equations for $L$ is a local section $\phi$ of $\pi\colon E\to M$ such that its first jet prolongation $j^1\phi$ satisfies

Classical Field Theory may be thought of a generalization of Classical Mechanics so as the dynamical system under study does not only depend on a one dimensional parameter, the time line, but on multi-dimensional one, space-time for instance. In this setting, a first order theory is described by a Lagrangian function, that is, a real function $L$ on the 1-jet bundle $J^1\pi$, where $\pi\colon E\to M$ is a fibration (Classical Mechanics corresponds to the particular case in which $M$ is the real line). A solution of the Euler-Lagrange equations for $L$ is a local section $\phi$ of $\pi\colon E\to M$ such that its first jet prolongation $j^1\phi$ satisfies
\[ \rbra\pp[L]{u^\alpha}-\dd{x^i}\pp[L]{u^\alpha_i}\rket \circ (j^1\phi) = 0 \,. \]
Here, $(x^i, u^\alpha)$ are local coordinates on $E$ which are adapted to the fibration $\pi$ and $(x^i, u^\alpha, u^\alpha_i)$  are the corresponding local coordinates on $J^1\pi$.

$J^1\pi$ is an affine bundle modeled on the vector bundle $V(J^1\pi)=\pi^*(T^*M)\tensor\Vert(\pi)$. The dual bundle $J^1\pi^\circ$ of $V(J^1\pi)$ is the so-called restricted multimomentum bundle and a section of the canonical projection $\mu\colon J^1\pi^\dag \to J^1\pi^\circ$ is said to be a Hamiltonian section, where $J^1\pi^\dag$ is the extended multimomentum bundle, that is, the affine dual to $J^1\pi$. A section $\tau$ of the canonical projection $\pi^\circ_1\colon J^1\pi^\circ\to M$ is a solution of the Hamilton-De Donder-Weyl equations for $h$ if
\[ \pp[u^\alpha ]{x^i} = \pp[H]{p_\alpha^i} \circ \tau \,, \quad
   \pp[p_\alpha^i]{x^i} = -\pp[H]{u^\alpha} \circ \tau \]
where $(x^i, u^\alpha, p^i_\alpha)$ are local coordinates on $J^1\pi^\circ$ and
\[ h(x^i, u^\alpha, p^i_\alpha) = (x^i, u^\alpha, -H(x^i, u^\alpha, p^i_\alpha), p^i_\alpha) \,. \]

A geometric formulation of these equations may be given using the canonical multisymplectic structure $\Omega$ of the extended multimomentum bundle $J^1\pi^\dag$ and the homogeneous real function $\bar H$ on $J^1\pi^\dag$ induced by the Hamiltonian section $h$. It is the extended Hamiltonian formalism for classical field theories of first order (see \cite{Cmps10,LnMnzRmn07}). Hamilton-De Donder-Weyl equations may also be obtained in an intrinsic form using the volume form $\eta$ on $M$ and the multisymplectic structure on $J^1\pi^\circ$ induced by $\Omega$ and $h$. It is the restricted Hamiltonian formalism (see, for instance, \cite{Cmps10,CrnnCrmpIbrt91,MnzRmn99,MnzRmn00a,MnzRmn00b,GchtMngtSrdn97,LnMrnMrr96,Rmn09,RmnMdstVlrn11}).

An extension of the previous Lagrangian and Hamiltonian formalism for classical field theories may be developed for the more general case when the base manifold is not, in general, orientable (see, for instance, \cite{Cmps10,GchtMngtSrdn97}).

Other different geometric approaches to classical field theories have been proposed by several authors using polysymplectic, $k$--symplectic or $k$--cosymplectic structures (see \cite{RmnMdstVlrn11} and references therein). However, these formulations only cover some special types of classical field theories.

Now, a natural question arises: Is it possible to develop a geometrical description of Lagrangian and Hamiltonian field theories in the Lagrangian submanifold setting? A positive answer to this question was given in \cite{LnLcmbRdrg91} for the particular case of $k$--symplectic (and $k$--cosymplectic) field theory (see also \cite{Grbw10, RmnMdstVlrn08}).

On the other hand, a Tulczyjew's triple for multisymplectic classical field theories was proposed in \cite{LnMrtSnt03}. The triple is adapted to the restricted Hamiltonian formalism and interesting ideas are involved in its construction. However, the intrinsic character of some constructions in \cite{LnMrtSnt03} is not proved  and, in addition, the Hamiltonian section is involved in the definition of the Hamiltonian side of the triple (this is an important difference if we compare with the original Tulczyjew's triple for Classical Mechanics).

Very recently, a Tulczyjew's triple for classical field theories of first order was proposed in \cite{Grbw11}. This triple is adapted to the restricted Hamiltonian formalism and the basic concepts of the variational calculus are used in order to construct it. Multisymplectic geometry is not used in this approach.

In this paper, we discuss a description of Lagrangian and Hamiltonian classical field theories in terms of Lagrangian submanifolds of premultisymplectic manifolds. For this purpose, we will construct a Tulczyjew's triple associated with a fibration $\pi\colon E\to M$. For the sake of simplicity, we consider the case where $M$ is an oriented manifold with a fixed volume form. However, in the appendix we extend the construction for the more general case where $M$ is not necessarily an oriented manifold. Our triple is adapted to the extended Hamiltonian formalism. The multisymplectic structure of the extended multimomentum bundle and the involution $\ex_\nabla$ of the iterated jet bundle associated with $\pi$ and a linear connection $\nabla$ on $M$ (see \cite{KlrMdgn91,Mdgn89}) play an important role in this construction. We remark that $\ex_\nabla$ is used in order to construct the Lagrangian side of the triple and that $\ex_\nabla$ depends on the linear connection $\nabla$. However, the structural maps of the triple do not depend on the linear connection. In this sense, the triple is canonical.

The paper is structured as follows. In Section 2, we present some basic notions on (pre)multisymplectic vector spaces and manifolds. The multisymplectic formulation of classical field theory is reviewed in Section 3. In Section 4, we discuss the Lagrangian and Hamiltonian side of our Tulczyjew's triple for classical field theories of first order. In Section 5, we present our conclusions and a description of future research directions. The paper ends with an appendix where we propose an extension of the Tulczyjew's triple for a fibration $\pi\colon E\to M$ such that $M$ needs not to be orientable and therefore a volume form on $M$ is no longer required.

% MULTISYMPLECTIC STRUCTURES ---------------------------------------------------
\section{Multisymplectic structures}
In the same way in which Classical Mechanics is modeled on symplectic geometry, one of the geometric approaches of Classical Field Theories is multisymplectic geometry (for the definition and properties of a multisymplectic structure, see \cite{CntrIbrtLn99}). In this section, we give some basic notions on multisymplectic vector spaces and manifolds.

\subsection{Multisymplectic vector spaces}
Along this section, $V$ will denote a real vector space of finite dimension and $W\subset V$ a vector subspace of the former.

\begin{definition} \label{def:multisymplectic.vector}
A ($k+1$)--form  $\Omega$ on $V$ is said to be \emph{multisymplectic} if it is non-degenerate, that is, if the linear map
\begin{eqnarray*}
\flat_\Omega\colon
  V &\To    & \Lambda^kV^*\\
  v &\Mapsto& \flat_\Omega(v):=i_v\Omega
\end{eqnarray*}
is injective. In such a case, the pair $(V,\Omega)$ is said to be a \emph{multisymplectic} vector space of order $k+1$.
\end{definition}

\begin{definition} \label{def:orthogonal.and.isotropy}
Let $(V,\Omega)$ be a multisymplectic vector space of order $k+1$. Given a vector subspace $W\subset V$ of $V$, we define the \emph{$l$--orthogonal complement of $W$}, with $1\leq l\leq k$, as the subspace of $V$
\[ W^{\perp, l} := \cbra v \in V \,:\, i_{v \wedge w_1 \wedge\ldots\wedge w_l}\Omega=0 \,, \forall w_1,\ldots, w_l\in W \cket \,. \]
Moreover, we say that $W$ is
\begin{itemize}
\item \emph{$l$--isotropic} if $W\subset W^{\perp,l}$;
\item \emph{$l$--coisotropic} if $W^{\perp,l} \subset W$;
\item \emph{$l$--Lagrangian} if $W=W^{\perp,l}$;
\item \emph{multisymplectic} if $W\cap W^{\perp,k}=\{0\}$.
\end{itemize}
\end{definition}

It easily seen that $W^{\perp,1}\subseteq\ldots\subseteq W^{\perp,k}$. Besides, if $W$ is a multisymplectic subspace of $V$ then, the pullback  $i^*\Omega$ of the multisymplectic form $\Omega$ by the canonical inclusion $i\colon W\hookrightarrow V$ is a multisymplectic form on $W$ of degree $k+1$.

A natural problem that arises in some particular dynamical systems is that the dynamical form under study fails to be non-degenerate, which is a case of further study.

\begin{definition}
A \emph{premultisymplectic} structure of order $k+1$ on a real vector space $V$ of finite dimension is a ($k+1$)--form $\Omega$ on $V$. The pair $(V,\Omega)$ is said to be a premultisymplectic vector space of order $k+1$.
\end{definition}

Let $(V,\Omega)$ be a premultisymplectic vector space of order $k+1$. Then, the quotient vector space $\tilde V:=V/\ker\flat_\Omega$ admits a natural multisymplectic structure $\widetilde\Omega$ of order $k+1$ which is characterized by the following condition
\[ \mu^*\widetilde\Omega = \Omega \,, \]
where $\mu\colon V \to\tilde V = V/\ker\flat_\Omega$ is the canonical projection.

\begin{definition}
Let $(V,\Omega)$ be a premultisymplectic vector space of order $k+1$ and $W$ be a subspace of $V$. Given $1\leq l \leq k$, then $W$ is said to be $l$--isotropic  (resp., $l$--coisotropic, $l$--Lagrangian, multisymplectic) if $\mu(W)$ is an $l$--isotropic (resp., $l$--coisotropic, $l$--Lagrangian, multisymplectic) subspace of the multisymplectic vector space $\tilde V = V/\ker\flat_\Omega$.
\end{definition}

Note that $W$ is $l$--isotropic (resp., $l$--coisotropic, $l$--Lagrangian, multisymplectic) for $(V,\Omega)$ if, and only if, the quotient space $W/(W \cap\ker\flat_\Omega)$ is an $l$--isotropic  (resp., $l$--coisotropic, $l$--Lagrangian, multisymplectic) subspace of the multisymplectic vector space $(\tilde V,\tilde\Omega)$.

\subsection{Multisymplectic manifolds}
Along this section, both $P$ and $E$ represent real smooth manifolds of finite dimension.

\begin{definition}
A \emph{premultisymplectic structure of order $k+1$} on a manifold $P$ is a closed ($k+1$)--form $\Omega$ on $P$. If furthermore $(T_xP,\Omega(x))$ is multisymplectic for each $x\in P$, then $\Omega$ is said to be a \emph{multisymplectic structure of order $k+1$} on $P$. The pair $(P, \Omega)$ is called a \emph{(pre)multisymplectic manifold of order $k+1$}.
\end{definition}

The canonical example of a multisymplectic manifold is the bundle of forms over a manifold $E$, that is, the manifold $P=\Lambda^kE$.

\begin{example}[\textbf{The bundle of forms}] \label{ex:multisymplectic:canonical}
Let $E$ be a smooth manifold of dimension $n$, $\Lambda^kE$ be the bundle of $k$--forms on $E$ and $\nu\colon \Lambda^kE\to E$ be the canonical projection. The \emph{Liouville or tautological form of order $k$} is the $k$--form $\Theta$ over $\Lambda^kE$ given by
\[ \Theta(\omega)(X_1,\ldots,X_k) := \omega((T_\omega\nu)(X_1),\ldots, (T_\omega\nu)(X_k)), \]
for any $\omega\in\Lambda^kE$ and any $X_1,\ldots,X_k\in T_\omega(\Lambda^kE)$. Then, the \emph{canonical multisymplectic ($k+1$)--form} is
\[ \Omega := -\diff\Theta \,. \]

If $(y^i)$ are local coordinates on $E$ and $(y^i, p_{i_1\ldots i_k})$, with $1\leq i_1<\ldots<i_k\leq n$, are the corresponding induced coordinates on $\Lambda^kE$, then
\begin{gather}
\Theta = \sum_{i_1<\ldots<i_k}p_{i_1\ldots i_k}\dy^{i_1}\wedge\ldots\wedge\dy^{i_k} \,,
\intertext{and}
\Omega = \sum_{i_1<\ldots<i_k}-\diff p_{i_1\ldots i_k}\wedge\dy^{i_1}\wedge\ldots\wedge\dy^{i_k} \,.
\end{gather}
From this expression, it is immediate to check that $\Omega$ is really multisymplectic. \hfill$\lhd$
\end{example}

\begin{example}[\textbf{The bundle of horizontal forms}] \label{ex:multisymplectic:horizontal}
Let $\pi\colon E\to M$ be a fibration, that is, $\pi$ is a surjective submersion. Assume that $\dim M=m$ and $\dim E=m+n$. Given $1\leq r\leq n$, we consider the vector subbundle $\Lambda^k_rE$ of $\Lambda^kE$ whose fiber at a point $u \in E$ is the set of $k$--forms at $u$ that are $r$--horizontal with respect to $\pi$, that is, the set
\[ (\Lambda^k_rE)_u = \{ \omega\in\Lambda^k_uE \,:\, i_{v_r}\ldots i_{v_1}\omega=0 \quad \forall v_1,\ldots,v_r\in\Vert_u(\pi) \} \,, \]
where $\Vert_u\pi=\ker(T_u\pi)$ is the space of tangent vectors at $u\in E$ that are vertical with respect to $\pi$.

We denote by $\nu_r$, $\Theta_r$ and $\Omega_r$ the restriction to $\Lambda^k_rE$ of $\nu$, $\Theta$ and $\Omega$ (given in the previous example). We have that $(\Lambda^k_rE,\Omega_r)$ is a multisymplectic manifold. The case in which $r=2$ and $k=m$ is of furthermost importance for multisymplectic field theory.

Let $(x^i,u^\alpha)$ denote adapted coordinates on $E$, then they induce coordinates $(x^i, u^\alpha, p, p^i_\alpha)$ on $\Lambda^m_2E$ such that any element $\omega\in\Lambda^m_2E$ has the form $\omega=p\dmx + p^i_\alpha\du^\alpha\wedge\dmxi$, where $\dmx = \dx^1\land ...\land \dx^m$ and $\dmxi = i_{\pp{x^i}}\dmx$. Therefore, we have that $\Theta_2$ and $\Omega_2$ are locally given by the expressions
\begin{gather}
\label{eq:form:tautological:coord}
\Theta_2 = p\dmx + p^i_\alpha\du^\alpha\wedge\dmxi \,,
\intertext{and}
\label{eq:form:mutlisymplectic:coord}
\Omega_2 = -\diff p\wedge\dmx - \diff p^i_\alpha\wedge\du^\alpha\wedge\dmxi \,.
\end{gather}
\hfill$\lhd$
\end{example}

Before ending this section, we recall the definition of $l$--Lagrangian submanifolds of a multisymplectic manifold and the natural extension of this notion to premultisymplectic manifolds.

\begin{definition} \label{def:lagrangian:submanifold}
Let $(P, \Omega)$ be a multisymplectic manifold of order $k+1$. A submanifold $S$ of $P$ is \emph{$l$--Lagrangian}, with $1\leq l\leq k$, if $T_xS$ is an $l$--Lagrangian subspace of the multisymplectic vector space $(T_xP, \Omega(x))$, for all $x\in S$. In other words, $S$ is an $l$--Lagrangian submanifold of $P$ if
\begin{equation}
(T_xS)^{\perp,l} = T_xS, \quad \forall x\in S \,.
\end{equation}
\end{definition}

\begin{definition} \label{def:lagrangian:submanifold:pre}
Let $(P, \Omega)$ be a premultisymplectic manifold of order $k+1$ and $\flat_\Omega\colon TP\to\Lambda^kT^*P$ be the corresponding vector bundle morphism. A submanifold $S$ of $P$ is \emph{$l$--Lagrangian}, with $1\leq l\leq k$, if $T_xS/(T_xS\cap\ker\flat_{\Omega(x)})$ is an $l$--Lagrangian subspace of the multisymplectic vector space $T_xP/\ker\flat_{\Omega(x)}$, for all $x \in S$.
\end{definition}

% MULTISYMPLECTIC FORMULATION OF CFT -------------------------------------------
\section{Multisymplectic formulation of Classical Field Theory} \label{sec:cft:multisymplectic}
In this section, we recall the basics of the geometric formulation of Classical Field Theory within a multisymplectic framework of jet bundles. The theory is set in a configuration fiber bundle $\pi\colon E\to M$, whose sections represent the fields of the system. Then, one may choose to develop a Lagrangian formalism, by considering a Lagrangian density $\LL\colon J^1\pi\to\Lambda^mM$ (or a Lagrangian function $L\colon J^1\pi\to\RR$, if a volume form on $M$ has been fixed), and derive the Euler-Lagrange equations. Or one may choose to develop a Hamiltonian formalism, by considering a Hamiltonian density $\HH\colon J^1\pi^\dag\to\Lambda^m M$  (or a Hamiltonian section $h\colon J^1\pi^\circ\to J^1\pi^\dag$), and derive the Hamilton's equations. The literature on this subject is vast. For a more concise study, the interested reader is referred, for instance, to \cite{Cmps10,CmpsLnMrtVnk09,CrnnCrmpIbrt91,LnMnzRmn07,MnzRmn00a,LnMrnMrr96,Rmn09,Sndr89}.

From here on, $\pi\colon E\to M$ will always denote a fiber bundle of rank $n$ over an $m$--dimensional manifold, \ie\ $\dim M=m$ and $\dim E=m+n$. Fibered coordinates on $E$ will be denoted by $(x^i, u^\alpha)$, $1\leq i \leq m$, $1 \leq \alpha \leq n$; where $(x^i)$ are local coordinates on $M$. The shorthand $\dmx = \dx^1\wedge\ldots\wedge\dx^m$ will represent the local volume form that $(x^i)$ defines, but we also use the notation $\dmxi = i_{\pp*{x^i}}\dx^1\wedge\ldots\wedge\dx^m$ for the contraction with the coordinate vector field. Many bundles will be considered over $M$ and $E$, but all of them vectorial or affine. For these bundles, we will only consider natural coordinates. In general, indexes denoted with lower case Latin letters (resp. Greek letters) will range between 1 and $m$ (resp. 1 and $n$). The Einstein sum convention on repeated crossed indexes is always understood.

Furthermore, we assume $M$ to be orientable with fixed orientation, together with a determined volume form $\eta$. Its pullback to any bundle over $M$ will still be denoted $\eta$, as for instance $\pi^*\eta$. In addition, local coordinates on $M$ will be chosen compatible with $\eta$, which means such that $\dmx=\eta$.

The most part of the objects that will be defined on the subsequent sections will depend in some way or another on the volume form $\eta$. An alternative definition independent of the volume form may be given for each one of them (as an overall volume independent theory). In particular, the aim of the Appendix section is to provide a construction of the Tulczyjew's triple for field theory independent of the volume form.

%%% LAGRANGIAN FORMALISM -------------------------------------------------------
\subsection{Lagrangian formalism} \label{sec:cft:lagrangian:formalism}
As already mentioned, the Lagrangian formulation of Classical Field Theory is stated on the \emph{first jet manifold} $J^1\pi$ of the configuration bundle $\pi\colon E\to M$. This manifold is defined as the collection of tangent maps of local sections of $\pi$. More precisely,
\[ J^1\pi := \cbra T_x\phi \,:\, \phi\in\Sections_x(\pi),\ x\in M \cket \,. \]
The elements of $J^1\pi$ are denoted $j^1_x\phi$ and called the \emph{1st-jet of $\phi$ at $x$}. Adapted coordinates $(x^i, u^\alpha)$ on $E$ induce coordinates $(x^i, u^\alpha, u^\alpha_i)$ on $J^1\pi$ such that $u^\alpha_i(j^1_x\phi)=\pp*[\phi^\alpha]{x^i}|_x$. It is clear that $J^1\pi$ fibers over $E$ and $M$ through the canonical projections $\pi_{1,0}\colon J^1\pi \to E$ and $\pi_1\colon J^1\pi \to M$, respectively. In local coordinates, these projections are given by $\pi_{1,0}(x^i, u^\alpha, u^\alpha_i) = (x^i, u^\alpha)$ and $ \pi_1(x^i, u^\alpha, u^\alpha_i) = (x^i)$.

Fiberwise, $J^1\pi$ may be seen as a set of linear maps, namely for each $u\in E$,
\[ J^1_u\pi = \{ z\in\Lin(T_{\pi(u)}M,T_uE) \,:\, T_u\pi \circ z = \id_{T_{\pi(u)}M} \} \,, \]
which is an affine space modeled on the vector space
\[ V_u(J^1\pi) = T_{\pi(u)}^*M \tensor\Vert_u(\pi)= \{ \bar z\in\Lin(T_{\pi(u)}M,T_uE) \,:\, T_u\pi \circ \bar z = 0 \} \,. \]
Thus, the first jet manifold $J^1\pi$ is an affine bundle over $E$ modeled on the vector bundle $V(J^1\pi) = \pi^*(T^*M)\tensor_E\Vert(\pi)$.

The dynamics of a Lagrangian field system are governed by a \emph{Lagrangian density}, a fibered map $\LL\colon J^1\pi\to\Lambda^mM$ over $M$. The real valued function $L\colon J^1\pi \to \RR$ that satisfies $\LL=L\eta$ is called the \emph{Lagrangian function}. Both Lagrangians permit to define the so-called \emph{Poincaré-Cartan forms}:
\begin{equation} \label{eq:pc-forms}
\Theta_\LL = L\eta + \abra S_\eta,\diff L\aket \in\Forms^m(J^1\pi) \quand \Omega_\LL= -\diff\Theta_\LL \in\Forms^{m+1}(J^1\pi) \,,
\end{equation}
where $S_\eta$ is a special and canonical structure of $J^1\pi$ called \emph{vertical endomorphism} and whose local expression is
\begin{equation} \label{eq:vertical:endomorphism:eta}
S_\eta= (\du^\alpha-u^\alpha_j\dx^j)\wedge\dmxi \tensor \pp{u^\alpha_i} \,.
\end{equation}
In local coordinates, the Poincaré-Cartan forms read
\begin{align}
% \label{eq:dl}
% \diff\LL =& \rbra \pp[L]{u^\alpha}\du^\alpha + \pp[L]{u^\alpha_i}\du^\alpha _i \rket \wedge \dmx \,, \\
\label{eq:pc-mform:coord}
\Theta_\LL =& \rbra L-u^\alpha_i\pp[L]{u^\alpha_i} \rket \wedge \dmx + \pp[L]{u^\alpha_i}\du^\alpha\wedge \dmxi \,, \\
\label{eq:pc-m1form:coord}
\Omega_\LL =& -(\du^\alpha-u^\alpha_j\dx^j) \wedge \rbra \pp[L]{u^\alpha}\dmx - \diff\rbra\pp[L]{u^\alpha_i}\rket\wedge \dmxi\rket \,.
\end{align}
It is important to note that, even though the Poincaré-Cartan $m$--form (and therefore, the Poincaré-Cartan ($m+1$)--form) depends on the Lagrangian function $L$ and the volume form $\eta$, actually it only depends on the Lagrangian density $\LL=L\eta$ (see Appendix \ref{sec:t3cft:novol}).

A \emph{critical point} of $\LL$ is a (local) section $\phi$ of $\pi$ such that
\[(j^1\phi)^*(i_X \Omega_\LL) = 0,\]
for any vector field $X$ on $J^1\pi$. A straightforward computation would show that this implies that
\begin{equation} \label{eq:euler-lagrange:coord}
(j^1\phi)^*\rbra \pp[L]{u^\alpha}-\dd{x^i}\pp[L]{u^\alpha_i}\rket = 0 \,,\ 1\leq \alpha\leq n \,.
\end{equation}
The above equations are called \emph{Euler-Lagrange equations}.

%%% HAMILTONIAN FORMALISM ------------------------------------------------------
\subsection{Hamiltonian formalism} \label{sec:cft:hamiltonian:formalism}
The dual formulation of the Lagrangian formalism is the Hamiltonian one, which is set in the affine dual bundles of $J^1\pi$. The \emph{(extended) affine dual bundle} $J^1\pi^\dag$ is the collection of real-valued affine maps defined on the fibers of $\pi_{1,0}:J^1\pi\to E$, namely
\[ (J^1\pi)^\dag := \Aff(J^1\pi,\RR) = \cbra A\in\Aff(J^1_u\pi,\RR) \,\colon\, u\in E \cket \,. \]
The \emph{(reduced) affine dual bundle} $J^1\pi^\circ$ is the quotient of $J^1\pi^\dag$ by constant affine maps, namely
\[ (J^1\pi)^\circ := \Aff(J^1\pi,\RR)/\{f\colon E\to\RR\} \,. \]
It is again clear that $J^1\pi^\dag$ and $J^1\pi^\circ$ are fiber bundles over $E$ but, in contrast to $J^1\pi$, they are vector bundles. Moreover, $J^1\pi^\dag$ is a principal $\RR$--bundle over $J^1\pi^\circ$. The different canonical projections are denoted $\pi_{1,0}^\dag\colon J^1\pi^\dag\to E$, $\pi_{1,0}^\circ\colon J^1\pi ^\circ\to E$ and $\mu\colon J^1\pi^\dag\to J^1\pi^\circ$. The natural pairing between the elements of $J^1\pi^\dag$ and those of $J^1\pi$ will be denoted by the usual angular bracket,
\[ \pairing \colon J^1\pi^\dag\times_EJ^1\pi \To \RR \,. \]
We note here that $J^1\pi^\circ$ is isomorphic to the dual bundle of $V(J^1\pi) = \pi^*(T^*M)\tensor_E\Vert(\pi)$.

Besides of defining the affine duals of $J^1\pi$, we also consider the \emph{extended \emph{and} reduced multimomentum spaces}
\[ \MM\pi := \Lambda^m_2E \quand \MM^\circ\pi := \Lambda^m_2E/\Lambda^m_1E \,. \]
By definition, these spaces are vector bundles over $E$ and we denote their canonical projections $\nu\colon\MM\pi\to E$, $\nu^\circ\colon\MM^\circ\pi\to E$ and $\mu\colon\MM\pi\to\MM^\circ\pi$ (some abuse of notation here). Again, $\mu\colon\MM\pi\to\MM^\circ\pi$ is a principal $\RR$--bundle. We recall from Example \ref{ex:multisymplectic:horizontal} that $\MM\pi$ has a canonical multisymplectic structure which we denote $\Omega$. On the contrary, $\MM^\circ\pi$ has no canonical multisymplectic structure, but $\Omega$ can still be pulled back by any section of $\mu\colon\MM\pi\to\MM^\circ\pi$ to give rise to a multisymplectic structure on $\MM^\circ\pi$.

An interesting and important fact is how the four bundles we have defined so far are related. We have that
\begin{equation} \label{eq:dual:identification}
J^1\pi^\dag \cong \MM\pi \quand J^1\pi^\circ \cong \MM^\circ\pi \,,
\end{equation}
although these isomorphisms depend on the base volume form $\eta$. In fact, the bundle isomorphism $\Psi:\MM\pi\to J^1\pi^\dag$ is characterized by the equation
\[ \abra\Psi(\omega),j^1_x\phi\aket\eta = \phi^*_x(\omega) \,,\quad \forall j^1_x\phi\in J^1_{\nu(\omega)}\pi \,,\quad \forall\omega\in \MM\pi \,. \]
We therefore identify $\MM\pi$ with $J^1\pi^\dag$ (and $\MM^\circ\pi$ with $J^1\pi^\circ$) and use this isomorphism to pullback the duality nature of $J^1\pi^\dag$ to $\MM\pi$.

Adapted coordinates in $\MM\pi$ (resp. $\MM^\circ\pi$) will be of the form $(x^i,u^\alpha,p,p^i_\alpha)$ (resp. $(x^i,u^\alpha,p^i_\alpha)$), such that
\[ p\dmx+p^i_\alpha\du^\alpha\wedge\dmxi\in\Lambda^m_2E \quad (p\dmx\in\Lambda^m_1E) \,. \]
Under these coordinates, the canonical projections have the expression
\[ \nu(x^i,u^\alpha,p,p^i_\alpha)=(x^i,u^\alpha) \,,\quad
   \nu^\circ(x^i,u^\alpha,p^i_\alpha)=(x^i,u^\alpha) \quand
   \mu(x^i,u^\alpha,p,p^i_\alpha)=(x^i,u^\alpha,p^i_\alpha) \,; \]
and the paring takes the form
\[ \abra(x^i,u^\alpha,p,p^i_\alpha),(x^i,u^\alpha,u^\alpha_i)\aket = p + p^i_\alpha u^\alpha \,. \]
We also recall the local description of the canonical multisymplectic form $\Omega$ of $\MM\pi$,
\[ \Omega = -\diff p\wedge\dmx - \diff p^i_\alpha\wedge\du^\alpha\wedge\dmxi \,. \]

Now, we focus on the principal $\RR$--bundle structure of $\mu\colon\MM\pi\to\MM^\circ\pi$. This structure arises from the $\RR$--action
\begin{align*}
                \RR\times\MM\pi &\,\To    \, \MM\pi\\
                     (t,\omega) &\,\Mapsto\, t\cdot\eta_{\nu(\omega)} + \omega \,.
\intertext{In coordinates,}
(t,(x^i,u^\alpha,p,p^i_\alpha)) &\,\Mapsto\, (x^i,u^\alpha,t+p,p^i_\alpha) \,.
\end{align*}
We will denote by $V_\mu\in\Fields(\MM\pi)$ the infinitesimal generator of the action of $\RR$ on $\MM\pi$, which in coordinates is nothing else but $V_\mu=\pp{p}$. Since the orbits of this action are the fiber of $\mu$, $V_\mu$ is also a generator of the vertical bundle $\Vert(\mu)$.

The dynamics of a Hamiltonian field system are governed by a \emph{Hamiltonian section}, a section $h\colon\MM^\circ\pi\to\MM\pi$ of $\mu\colon\MM\pi\to\MM^\circ$. In presence of the base volume form $\eta$, the set of Hamiltonian sections $\Sections(\mu)$ is in one-to-one correspondence with the set of functions $\{\bar H\in\Smooth(\MM\pi) : V_\mu(\bar H)=1 \}$ and with the set of \emph{Hamiltonian densities}, fibered maps $\HH\colon\MM\pi\to\Lambda^mM$ over $M$ such that $i_{V_\mu}\diff\HH=\eta$. Given a Hamiltonian section $h\colon\MM^\circ\pi\to\MM\pi$, the corresponding Hamiltonian density is
\[ \HH(\omega) = \omega - h(\mu(\omega)) \,,\ \forall\omega\in\MM\pi \,. \]
Conversely, given a Hamiltonian density $\HH\colon\MM\pi\to\Lambda^mM$, the corresponding Hamiltonian section is characterized by the condition
\[ \mathrm{im}h=\HH^{-1}(0) \,. \]
Obviously, $\HH=\bar H\eta$. In adapted coordinates,
\begin{align}
\label{eq:hamiltonian:section:coord}
h(x^i, u^\alpha, p^i_\alpha) &= (x^i, u^\alpha, p = -H(x^i, u^\alpha, p^i_\alpha), p^i_\alpha) \,,\\
\label{eq:hamiltonian:density:coord}
\HH(x^i, u^\alpha, p, p^i_\alpha) &= \rbra p + H(x^i, u^\alpha, p^i_\alpha)\rket\dmx \,.
\end{align}
The locally defined function $H$ is called the \emph{Hamiltonian function} and it must not be confused with the globally defined function $\bar H$ such that $\HH=\bar H\eta$.

A \emph{critical point} of $\HH$ is a (local) section $\tau$ of $\pi\circ\nu\colon\MM\pi\to M$ that satisfies the \emph{(extended) Hamilton-De Donder-Weyl equation}
\begin{equation} \label{eq:hamilton-dedonder-weyl}
\tau^*i_X(\Omega+\diff\HH) = 0 \,,
\end{equation}
for any vector field $X$ on $\MM\pi$.
%where $X_\eta$ is the section of the vector bundle $\Lambda^mTM\to M$ characterized by $\eta(X_\eta)=1$.
A \emph{critical point} of $h$ is a (local) section $\tau$ of $\pi\circ\nu^\circ\colon\MM^\circ\pi\to M$ that satisfies the \emph{(reduced) Hamilton-De Donder-Weyl equation}
\begin{equation} \label{eq:hamilton-dedonder-weyl}
(h\circ\tau)^*i_X(\Omega+\diff\HH) = 0 \,,
\end{equation}
for any vector field $X$ on $\MM\pi$.
A straightforward but cumbersome computation shows that both are equivalent to the following set of local equations known as \emph{Hamilton's equations}:
\begin{equation} \label{eq:hamilton:coord}
\pp[(u^\alpha \circ \tau)]{x^i} = \pp[H]{p_\alpha^i} \circ \tau \,, \qquad
\pp[(p_\alpha^i\circ \tau)]{x^i} = -\pp[H]{u^\alpha} \circ \tau \,. %\quad
%\pp[(p \circ \tau)]{x^i} = -\pp[(H \circ \tau)]{x^i}.
\end{equation}
%Note that, from the last group of equations, we deduce that $p\circ\tau = -H\circ\tau + \textrm{cst.}$

%%% EQUIVALENCE BETWEEN BOTH FORMALISMS ----------------------------------------
\subsection{Equivalence between both formalisms} \label{sec:cft:equivalence}
In this section, we present the equivalence between the Lagrangian and Hamiltonian formalisms of classical field theories for the case when the Lagrangian function is (hyper)regular (see below).

Let $\LL$ be a Lagrangian density. The \emph{(extended) Legendre transform} is the bundle morphism $\Leg_\LL\colon J^1\pi\to\MM\pi$ over $E$ defined as follows:
\begin{equation} \label{eq:legendre.transform.ext}
\Leg_\LL(j^1_x\phi)(X_1,\ldots,X_m) := (\Theta_\LL)_{j^1_x\phi}(\widetilde{X}_1,\ldots,\widetilde{X}_m),
\end{equation}
for all $j^1_x\phi \in J^1\pi $ and $X_i \in T_{\phi(x)}E$, where $\widetilde{X}_i \in T_{j^1_x\phi}J^1\pi$ are such that $T\pi_{1,0}(\widetilde{X}_i)=X_i$. The \emph{(reduced) Legendre transform} is the composition of $\Leg_\LL$ with $\mu$, that is, the bundle morphism
\begin{equation} \label{eq:legendre.transform.red}
\leg_\LL := \mu \circ \Leg_\LL\colon J^1\pi \to \MM^\circ\pi \,.
\end{equation}
In local coordinates,
\begin{eqnarray}
\label{eq:legendre:transform:ext:coord}
\Leg_\LL(x^i, u^\alpha, u^\alpha_i) &=& \rbra x^i, u^\alpha, L-\pp[L]{u^\alpha_i}u^\alpha_i, \pp[L]{u^\alpha_i}\rket \,, \\
\label{eq:legendre:transform:red:coord}
\leg_\LL(x^i, u^\alpha, u^\alpha_i) &=& \rbra x^i, u^\alpha, \pp[L]{u^\alpha_i}\rket \,,
\end{eqnarray}
where $L$ is the Lagrangian function associated to $\LL$, \ie\ $\LL=L\eta$.

From the definitions, we deduce that $(\Leg_\LL)^*(\Theta)= \Theta_\LL, (\Leg_\LL)^*(\Omega) = \Omega_\LL$, where $\Theta$ is the Liouville $m$--form on $\MM\pi$ and $\Omega$ is the canonical multisymplectic ($m+1$)--form. In addition, we have that the Legendre transformation $\leg_\LL\colon J^1\pi\to \MM^\circ\pi$ is a local diffeomorphism, if and only if, the Lagrangian function $L$ is \emph{regular}, that is, the Hessian $\rbra\frac{\partial^2L}{\partial u^\alpha_i\partial u^\beta_j}\rket$ is a regular matrix. When $\leg_\LL\colon J^1\pi\to \MM^\circ\pi$ is a global diffeomorphism, we say that the Lagrangian $L$ is \emph{hyper-regular}. In this case, we may define the Hamiltonian section $h\colon \MM^\circ\pi\To \MM\pi$
\begin{equation} \label{eq:legendre:hamiltonian:section}
 h = \Leg_\LL\circ\leg_\LL^{-1} \,,
\end{equation}
whose associated Hamiltonian density is
\begin{equation} \label{eq:legendre:hamiltonian:density}
 \HH(\omega) = \abra\omega,\leg_\LL^{-1}(\mu(\omega))\aket\eta - (\LL\circ\leg_\LL^{-1})(\mu(\omega)) \,,\quad \forall\omega\in\MM\pi \,.
\end{equation}
In coordinates,
\begin{equation} \label{eq:legendre:hamiltonian:section:coord}
h(x^i, u^\alpha, p^i_\alpha) = (x^i, u^\alpha, L(x^i,u^\alpha,u^\alpha_i) - p^i_\alpha u^\alpha_i, p^i_\alpha) \,,
\end{equation}
where $u^\alpha_i=u^\alpha_i(\leg_\LL^{-1}(x^i, u^\alpha, p^i_\alpha))$. Accordingly,
\begin{equation} \label{eq:legendre:hamiltonian:density:coord}
\HH(x^i, u^\alpha, p, p^i_\alpha) = \rbra p + p^i_\alpha u^\alpha_i - L \rket\dmx
\end{equation}
and the Hamiltonian function is
\begin{equation} \label{eq:legendre:hamiltonian:function:coord}
H(x^i, u^\alpha, p^i_\alpha) = p^i_\alpha u^\alpha_i - L \,.
\end{equation}

\begin{theorem}[\refer\ \cite{MnzRmn00a,LnMrnMrr96,Rmn09}]
Assume $\LL$ is a hyper-regular Lagrangian density. If $\phi$ is a solution of the Euler-Lagrange equations for $\LL$, then $\omega=\leg_\LL\circ j^1\phi$ is a solution of the Hamilton's equations for $h$. Conversely, if $\omega$ is a solution of the Hamilton's equations for $h$, then $\leg_\LL^{-1}\circ\,\omega$ is of the form $j^1\phi$, where $\phi$ is a solution of the Euler-Lagrange equations for $\LL$.
\end{theorem}

% TULCZYJEW'S TRIPLE -----------------------------------------------------------
\section{Tulczyjew's triple} \label{sec:t3cft}
In 1976, W. Tulczyjew presented in a couple of papers \cite{Tlcz76a,Tlcz76b} a complete geometric construction relating the triple $T(T^*Q)$, $T^*(TQ)$ and $T^*(T^*Q)$. The canonical symplectic structures of these spaces are carried between them maintaining their symplectic character. However, the relevance of his work is not only this but the fact that, when one introduces dynamics into the picture, a Lagrangian and a Hamiltonian system, then the dynamics are lifted to this triple giving the natural correspondence between these two formalisms and showing again the very geometric nature of Lagrangian and Hamiltonian mechanics.

In this section, we will obtain the main results of the paper. In fact, we will introduce a Tulczyjew's triple for Classical Field Theory associated with a fibration $\pi\colon E \to M$, where $M$ is an oriented manifold with a fixed volume form $\eta$. First, we will present the Lagrangian side of the triple and then the Hamiltonian side. Finally, we will introduce the equivalence between the two formalisms in this setting.

%%% LAGRANGIAN FORMALISM -------------------------------------------------------
\subsection{Lagrangian formalism} \label{sec:t3cft:lagrangian}
Recall that if $Q$ is a smooth manifold, there is a canonical vector bundle isomorphism $A_Q\colon T(T^*Q)\to T^*(TQ)$, named after Tulczyjew's work \cite{Tlcz76b}, whose local representation is
\[ A_Q(q^i, p_i, \dot q^i, \dot p_i) = (q^i, \dot q^i, \dot p_i, p_i) \,, \]
where $(q^i)$ are local coordinates on $Q$ and $(q^i, p_i)$ (respectively, $(q^i, p_i, \dot q^i, \dot p_i))$ are the corresponding adapted local coordinates on $T^*Q$ (respectively, $T(T^*Q)$). The Tulczyjew's isomorphism $A_Q$ is built by means of the canonical involution of $TTQ$, $\kappa_Q\colon TTQ\to TTQ$, and the tangent lift of the natural pairing $\pairing\colon T^*Q {\,}_{\pi_Q}\!\!\times_{\tau_Q} TQ\to\RR$, where $\tau_Q\colon TQ\to Q$ and $\pi_Q\colon T^*Q\to Q$ are the canonical tangent and cotangent projections.

Note that, while the second argument of the tangent pairing (the tangent lift of $\pairing$)
\[ T\pairing\colon TT^*Q {\,}_{T\pi_Q}\!\times_{T\tau_Q} TTQ\to\RR \]
fibers over $TQ$ through $T\tau_Q$, thanks to the involution $\kappa_Q$, the second argument of the nondegenerate pairing
\[ T\pairing \circ (id,\kappa_Q)\colon TT^*Q {\,}_{T\pi_Q}\!\!\times_{\tau_{TQ}} TTQ \to\RR \]
fibers over $TQ$ through $\tau_{TQ}$. In fact, $A_Q$ is the vector bundle isomorphism that this later map induces from $T(T^*Q)$ to $T^*(TQ)$.

In what follows, we will discuss the construction of the Tulczyjew's morphism for jet bundles. We shall define an affine morphism (between the proper spaces) by mimicking Tulczyjew's construction. Therefore, we will need an involution of the iterated jet bundle $J^1\pi_1$ and a ``tangent'' pairing.

Like the iterated tangent bundle $TTQ$ of a manifold $Q$, the iterated jet bundle $J^1\pi_1$ of a fiber bundle $\pi\colon E\to M$ has two different affine structures over $J^1\pi$. The first is the one in which we think of $J^1\pi_1$ as the 1st-jet bundle of the fiber bundle $\pi_1\colon J^1\pi\to M$, that is the affine bundle $(\pi_1)_{1,0}\colon J^1\pi_1\to J^1\pi$ modeled over the vector bundle $V(J^1\pi_1)=\pi^*_1(T^*M)\tensor_{J^1\pi}\Vert(\pi_1)$. The second is the one in which we think of $J^1\pi_1$ as the 1st-jet prolongation of the morphism $\pi_{1,0}\colon J^1\pi\to E$, that is the affine bundle $j^1(\pi_{1,0})\colon J^1\pi_1\to J^1\pi$ modeled over the vector bundle $J^1(V(J^1\pi))$.

\begin{figure}[h!]
\begin{minipage}{0.5\textwidth}
\[ \xymatrix{
     J^1\pi_1 \ar[rr]^{j^1(\pi_{1,0})} \ar[dd]_{(\pi_1)_{1,0}} && J^1\pi \ar[dd]^{\pi_{1,0}} \\ \\
     J^1\pi \ar[rr]^{\pi_{1,0}} && E
} \]
\caption{The iterated jet bundle} \label{fig:iterated.jet.bundle}
\end{minipage}
\begin{minipage}{0.45\textwidth}
\[ \xymatrix{
     J^1\pi_1 \ar[ddr]_{(\pi_1)_{1,0}} \ar@{<->}[rr]^{\ex_\nabla}  && J^1\pi_1 \ar[ddl]^{j^1(\pi_{1,0})} \\ \\
     & J^1\pi &
} \]
\caption{The exchange map} \label{fig:exchange.map}
\end{minipage}
\end{figure}

But unlike the iterated tangent bundle $TTQ$, there isn't a canonical involution of the iterated jet bundle $J^1\pi_1$ of a fiber bundle $\pi\colon E\to M$ (besides of the identity map). However, Kol\'a\v{r} and Modugno showed in \cite{KlrMdgn91} that the natural involutions of the iterated jet bundle $J^1\pi_1$ depend on the symmetric linear connections of the base manifold $M$. In fact, given a symmetric linear connection $\nabla$ on $M$, they introduced an affine bundle isomorphism $\ex_\nabla$ over the identity of $J^1\pi$ from $(\pi_1)_{1,0}\colon J^1\pi_1\to J^1\pi$ to $j^1(\pi_{1,0})\colon J^1\pi_1\to J^1\pi$. In local coordinates, the \emph{exchange map} $\ex_\nabla:J^1\pi_1\to J^1\pi_1$ has the expression
\begin{equation} \label{eq:exchange.map.coord}
\ex_\nabla(x^i, u^\alpha, u^\alpha_i, \bar u^\alpha_j, u^\alpha_{ij}) = (x^i, u^\alpha, \bar u^\alpha_i, u^\alpha_j, u^\alpha_{ji} + (\bar u^\alpha_k-u^\alpha_k)\chr^k_{ji}) \,,
\end{equation}

where $\chr^k_{ij}$ are the Christoffel's symbols of the symmetric linear connection $\nabla$. Note that $\ex_\nabla \circ \ex_\nabla = id_{J^1\pi}$.

The definition of the ``tangent'' pairing needs of an extra element that also depends on the linear connection $\nabla$ and on the volume form $\eta$ too. We define the $\RR$--linear map $\diff^{\nabla, \eta}\colon \Smooth(M) \to \Forms^1(M)$ characterized by the condition
\begin{equation} \label{eq:dnablaeta}
\diff^{\nabla,\eta} f\tensor\eta = \nabla(f\cdot\eta) \,, \quad \forall f\in\Smooth(M) \,.
\end{equation}
For a coordinate system $(x^i)$ compatible with the volume form, \ie\ $\eta = \dmx$, we have that
\[ \diff^{\nabla,\eta} f = \rbra\pp[f]{x^i}-f\cdot\chr^j_{ij}\rket\dx^i \,. \]

Now, we consider the natural pairing between $\MM\pi$ and $J^1\pi$
\[ \begin{array}{rcl}
     \pairing\colon
       \MM\pi {\,}_\nu\!\!\times_{\pi_{1,0}} J^1\pi &\To    & \RR \\
       (\omega,j^1_x\phi) &\Mapsto& \abra\omega,j^1_x\phi\aket \quad\textrm{such that}\quad \phi^*_x(\omega)=\abra\omega,j^1_x\phi\aket\eta(x)
\end{array} \]
and we lift it to the map
\[ \begin{array}{rcl}
     \diff^{\nabla,\eta}\pairing\colon
       J^1(\pi \circ \nu) {\,}_{j^1\nu}\!\!\times_{j^1(\pi_{1,0})} J^1\pi_1 &\To    & T^*M \\
       (j^1_x\omega,j^1_x\sigma) &\Mapsto& \diff^{\nabla,\eta}\abra\omega,\sigma\aket(x) \, .
\end{array} \]
In local coordinates,
\[ \abra(x^i, u^\alpha, p, p_\alpha^i),\, (x^i, u^\alpha, u^\alpha_i)\aket = p+p_\alpha^iu^\alpha_i \, . \]
Therefore,
\begin{multline*}
\diff^{\nabla,\eta}\abra(x^i, u^\alpha, p, p_\alpha^i, \bar u^\alpha_j, p_j, p_{\alpha j}^i),\,(x^i, u^\alpha, u^\alpha_i, \bar u^\alpha_j, u^\alpha_{ij})\aket =\\
= \rbra p_k + p^i_{\alpha k}u^\alpha_i+p^i_\alpha u^\alpha_{ik} - (p+p^i_\alpha u^\alpha_i)\chr^j_{kj}\rket\dx^k \,.
\end{multline*}

At this point, we remark that the second argument of the prolonged pairing is the iterated jet bundle $J^1\pi_1$ fibering over $J^1\pi$ by $j^1(\pi_{1,0})$. Thus, by composition with the exchange map $\ex_\nabla$, we obtain a morphism
\[ \diff^{\nabla,\eta}\pairing\circ(\id\times\ex_\nabla)\colon J^1(\pi \circ \nu) {\,}_{j^1\nu}\!\!\times_{(\pi_1)_{1,0}} J^1\pi_1 \To  T^*M \]
where now the second argument is the iterated jet bundle $J^1\pi_1$ fibering over $J^1\pi$ by $(\pi_1)_{1,0}$. In local coordinates,
\begin{multline*}
\diff^{\nabla,\eta}\abra(x^i, u^\alpha, p, p_\alpha^i, u^\alpha_j, p_j, p_{\alpha j}^i), \ex_\nabla(x^i, u^\alpha, u^\alpha_i, \bar u^\alpha_j, u^\alpha_{ij})\aket =\\
= \rbra p_k + p^i_{\alpha k}\bar u^\alpha_i+p^i_\alpha(u^\alpha_{ki}+ (\bar u^\alpha_l-u^\alpha_l)\chr^l_{ki}) - (p+p^i_\alpha\bar u^\alpha_i)\chr^j_{kj}\rket\dx^k \,.
\end{multline*}
Besides, observe that $\diff^{\nabla,\eta}\pairing\circ(\id\times\ex_\nabla)$ is a vector valued map that takes values from a vector bundle, $J^1(\pi \circ \nu)$, and an affine bundle, $J^1\pi_1$. Hence it induces a vector bundle morphism over the identity of $J^1\pi$
\[ \widetilde{A^{\nabla,\eta}_\pi}\colon J^1(\pi\circ\nu) \To \Aff_{\pi_1}(J^1\pi_1,T^*M) \,, \]
where
\[ \Aff_{\pi_1}(J^1\pi_1, T^*M) = \cup_{z\in J^1\pi} \, \Aff(J^1_z\pi_1, T^*_{\pi_1(z)}M) \,. \]
For adapted coordinates $(x^i, u^\alpha, u^\alpha_i, \bar p_k, \bar p^j_{\alpha k}, \bar p^{ij}_{\alpha k})$ on the vector bundle $\Aff_{\pi_1}(J^1\pi_1, T^*M)$,
\begin{multline*}
\widetilde{A_\pi^{\nabla,\eta}}(x^i, u^\alpha, p, p_\alpha^i, u^\alpha_j, p_j, p_{\alpha j }^i) =\\
= (x^i, u^\alpha, u^\alpha_i,
\bar p_k = p_k - p^i_\alpha u^\alpha_l\chr^l_{ki} - p\chr^j_{kj},
\bar p^i_{\alpha k} = p^i_{\alpha k} + p^l_\alpha\chr^i_{kl} - p^i_\alpha\chr^j_{kj},
\bar p^{ij}_{\alpha k} = p^j_\alpha\delta^i_k) \,.
\end{multline*}

\begin{lemma} \label{th:affine.map.lemma}
Given an arbitrary fiber bundle $\pi_{E,M}\colon E\to M$, let $\pi_{F,E}\colon F\to E$ and $\pi_{V,M}\colon V\to M $ be an affine and a vector bundle, respectively. We have that, the vector bundles
\[ \Aff_{\pi_{E,M}}(F,V) := \cup_{y\in E}\, \Aff(F_y, V_{\pi_{E,M}(y)}) \,,\]
and
\[ F^\dag\tensor_E\pi_{E,M}^*(V) = \Aff(F,\RR)\tensor_E\pi_{E,M}^*(V) \]
are isomorphic.
\end{lemma}

\begin{proof}
The lemma's assertion is a simple algebraic fact. Given a point $x\in M$, lets fix a point $y\in E_x$ in its fiber, $\pi_{E,M}(y)=x$. Then, we define a map from $F^\dag_y \times V_x$ to $\Aff(F_y, V_x)$ as follows: Given $\omega\in F^\dag_y$ and $v\in V_x$, we consider the affine map $(\omega:v)\colon F_y\to V_x$ given by $(\omega:v)(z) = \omega(z)v$, for $z \in F_y$. The map
\[ (\omega,v)\in F^\dag_y \times V_x \Mapsto (\omega:v)\in\Aff(F_y, V_x) \]
is obviously bilinear and, therefore, it induces a linear map from the tensor product $F^\dag_y \tensor V_x$ to $\Aff(F_y, V_x)$. It only rests to prove that, in fact, it is an isomorphism.

Consider a basis $(v_j)$ of $V_x$ and a basis $\{1, \omega^i\}$ of $F_y^\dag$ dual to a reference system $(o, f_i)$ of $F_y$, \ie
\[ 1(o+r^if_i) = 1 \quand \omega^j(o+ r^if_i)= r^j \,,\ \forall r^i\in \RR \,. \]
Then $\{1\tensor v_j,\ \omega^i\tensor v_j\}$ is a basis of $F_y^\dag\tensor V_x$ and, as it is easy to check, its image $\{(1:v_j),\ (\omega^i:v_j)\}$ is a basis of $\Aff(F_y,V_x)$.
\end{proof}

Using the previous lemma, we transform the range space of $\widetilde{A_\pi^{\nabla,\eta}}$, obtaining a new morphism
\[ \widetilde{A_\pi^{\nabla,\eta}}\colon J^1(\pi \circ \nu)\To(J^1\pi_1)^\dag \tensor_{J^1\pi}(\pi_1)^*(T^*M) \,, \]
which we continue to denote by $ \widetilde{A_\pi^{\nabla,\eta}}$. Note that $(J^1\pi_1)^\dag$ may be identified with the bundle of $m$--forms $\MM\pi_1 = \Lambda^m_2J^1\pi$. Under this identification, the local expression of $\widetilde{A_\pi^{\nabla,\eta}}$ is
\begin{multline*}
\widetilde{A_\pi^{\nabla,\eta}}(x^i, u^\alpha, p, p_\alpha^i, u^\alpha_j, p_j, p_{\alpha j}^i)=\\
= (\bar p_k\dmx +\bar p^i_{\alpha k}\du^\alpha\wedge \dmxi + \bar p^{ij}_{\alpha k}\du^\alpha_i\wedge \dmxj) \tensor\dx^k \,,
\end{multline*}
where $\bar p_k = p_k - p^i_\alpha u^\alpha_l\chr^l_{ki} - p\chr^j_{kj}$, $\bar p^i_{\alpha k} = p^i_{\alpha k} + p^l_\alpha\chr^i_{kl} - p^i_\alpha\chr^j_{kj}$, $\bar p^{ij}_{\alpha k} = p^j_\alpha\delta^i_k$.

Considering the natural morphism $(J^1\pi_1)^\dag\tensor_{J^1\pi}(\pi_1)^*(T^*M)\to\Lambda^{m+1}_2J^1\pi$ given by the wedge product, we finally obtain the vector bundle morphism
\begin{equation}\label{Apieta}
\begin{array}{rcl}
                       A_\pi^\eta\colon J^1(\pi\circ\nu) &\To    & \Lambda^{m+1}_2J^1\pi\\
(x^i, u^\alpha, p, p_\alpha^i, u^\alpha_j, p_j, p_{\alpha j}^i) &\Mapsto& (p_{\alpha i}^i\du^\alpha + p_\alpha^i\du^\alpha_i) \wedge \dmx \,,
\end{array}
\end{equation}
which we call \emph{the Tulczyjew's morphism}. It is important to note that this morphism does no longer depend on the linear connection $\nabla$, while it still depends on the volume form $\eta$ even though it is not explicitly noted.

\begin{theorem} \label{th:tulczyjew:morphism:left}
The Tulczyjew's morphism $A_\pi^\eta\colon J^1(\pi\circ\nu) \to \Lambda^{m+1}_2J^1\pi$, locally given by Equation \eqref{Apieta}, is a vector bundle epimorphism over the identity of $J^1\pi$. Moreover, it is canonical in the sense that it only depends on the original fiber bundle $\pi\colon E\to M$ and, a pripori, on the base volume form $\eta$.
\end{theorem}

According to Example \ref{ex:multisymplectic:horizontal}, $\Lambda^{m+1}_2J^1\pi$ has a canonical multisymplectic structure given by the form
\begin{equation}
\Omega_{\Lambda^{m+1}_2J^1\pi } = -\diff\bar p_\alpha \wedge\du^\alpha\wedge\dmx -\diff\bar p_\alpha^i\wedge\du^\alpha_i\wedge\dmx \,,
\end{equation}
where $(x^i, u^\alpha, u^\alpha_i, \bar p_\alpha, \bar p_\alpha^i)$ are natural coordinates on $\Lambda^{m+1}_2J^1\pi$. This structure is pulled back to $J^1(\pi\circ\nu)$ by $A_\pi^\eta$ defining a premultisymplectic structure given by the ($m+2$)--form
\begin{equation} \label{eq:omega:tilde:coord}
\widetilde\Omega:=(A_\pi^\eta)^*(\Omega_{\Lambda^{m+1}_2J^1\pi}) = -\diff p^i_{\alpha i}\wedge\du^\alpha\wedge\dmx - \diff p^i_\alpha\wedge\du^\alpha_i\wedge\dmx \,.
\end{equation}
The local basis of the kernel of $\widetilde\Omega$ is generated by
\begin{equation} \label{eq:omega:tilde:kernel}
\ker\widetilde\Omega = \abra\cbra \pp{p}, \pp{p_j}, \pp{p^i_{\alpha j}} - \delta^i_j\pp{p^k_{\alpha k}} \cket\aket \,,
\end{equation}
where $1\leq k\leq m$ is a fixed index (no Einstein convention here, take for instance $k=1$).

Next, let $\LL\colon J^1\pi\to\Lambda^mM$ be an arbitrary Lagrangian density such that $\LL=L\eta$. Then, we obtain the following result.

\begin{proposition} \label{th:SL:lagrangian:manifold}
$S_\LL = (A_\pi^\eta)^{-1}(\diff\LL(J^1\pi))$ is an ($m+1$)--Lagrangian submanifold of the premultisymplectic manifold $(J^1(\pi \circ \nu), \widetilde\Omega)$.
\end{proposition}

\begin{proof}
From \eqref{Apieta}, we obtain that the submanifold $S_\LL$ is locally given by
\begin{equation} \label{eq:SL}
S_\LL= \cbra (x^i, u^\alpha, p, p^i_\alpha, u^\alpha_j, p_j, p^i_{\alpha j}) \,:\, p^i_\alpha=\pp[L]{u^\alpha_i}, \quad p^i_{\alpha i}=\pp[L]{u^\alpha} \cket \,,
\end{equation}
thus
\begin{equation} \label{eq:TSL}
TS_\LL = \abra\cbra X_i, U_\alpha, U_\alpha^i, \pp{p}, \pp{p_j}, \pp{p^i_{\alpha j}} - \delta^i_j\pp{p^1_{\alpha 1}} \cket\aket \,,
\end{equation}
where
\begin{align}
\label{eq:TSL:Xi}
X_i =& \pp{x^i}
+ \frac{\partial^2L}{\partial x^i \partial u^\beta_j}\pp{p^j_\beta}
+ \frac{\partial^2L}{\partial x^i\partial u^\beta}\pp{p^1_{\beta 1}} \,, \\
\label{eq:TSL:Ua}
U_\alpha =& \pp{u^\alpha}
+ \frac{\partial^2L}{\partial u^\alpha\partial u^\beta_j}\pp{p^j_\beta}
+ \frac{\partial^2L}{\partial u^\alpha\partial u^\beta}\pp{p^1_{\beta 1}} \,, \\
\label{eq:TSL:Uai}
U_\alpha^i =& \pp{u^\alpha_i}
+ \frac{\partial^2L}{\partial u^\alpha_i\partial u^\beta_j}\pp{p^j_\beta}
+ \frac{\partial^2L}{\partial u^\alpha_i\partial u^\beta}\pp{p^1_{\beta 1}} \,.
\end{align}
It follows from Equation \eqref{eq:omega:tilde:kernel} that
\[ \ker\widetilde\Omega(j^1_x\omega) \subseteq T_{j^1_x\omega}S_\LL \,,\quad \forall j^1_x\omega \in S_\LL \,; \]
and, from Definitions \ref{def:lagrangian:submanifold} and \ref{def:lagrangian:submanifold:pre}, we must prove that
\[ \frac{T_{j^1_x\omega}S_\LL}{\ker\widetilde\Omega(j^1_x\omega)} = \rbra  \frac{T_{j^1_x\omega}S_\LL}{\ker\widetilde\Omega(j^1_x\omega)} \rket^{\perp,m+1} \,,\quad \forall j^1_x\omega \in S_\LL \,. \]

If $i_{S_\LL}\colon S_\LL\to J^1(\pi\circ\nu)$ is the natural inclusion then, using \eqref{eq:omega:tilde:coord} and \eqref{eq:SL}, we deduce that
\[i_{S_\LL}^*\widetilde\Omega = -d \rbra\pp[L]{u^\alpha} \du^\alpha + \pp[L]{u^\alpha_i}\du^\alpha_i\rket \land \dmx = \diff\rbra \pp[L]{x^i}\dx^i\rket \land \dmx = 0\]
which implies that
\[ \frac{T_{j^1_x\omega}S_\LL}{\ker\widetilde\Omega(j^1_x\omega)} \subseteq \rbra  \frac{T_{j^1_x\omega}S_\LL}{\ker\widetilde\Omega(j^1_x\omega)} \rket^{\perp,m+1} \,. \]

Instead of showing the converse inclusion, we will show that there is no gap in between. For this, we first note that
\[ \frac{T_{j^1_x\omega}S_\LL}{\ker\widetilde\Omega(j^1_x\omega)} = \abra\cbra X_i, U_\alpha, U_\alpha^i \cket\aket \,, \]
where $X_i$, $U_\alpha$, and $U_\alpha^i$ are given by Equations \eqref{eq:TSL:Xi}, \eqref{eq:TSL:Ua} and \eqref{eq:TSL:Uai}. Besides,
\[ T_{j^1_x\omega}J^1(\pi\circ\nu) = \abra\cbra X_i, U_\alpha, U_\alpha^i, \pp{p}, \pp{p^i_\alpha}, \pp{p_j}, \pp{p^i_{\alpha j}} \cket\aket \,, \]
therefore
\[ \frac{T_{j^1_x\omega}J^1(\pi\circ\nu)}{\ker\widetilde\Omega(j^1_x\omega)} = \abra\cbra X_i, U_\alpha, U_\alpha^i, \pp{p^i_\alpha}, \pp{p^i_{\alpha i}} \cket\aket \,, \]
where $i$ goes from 1 to $m$ and the Einstein notation is not considered.

Now, if $\rbra T_{j^1_x\omega}S_\LL / \ker\widetilde\Omega(j^1_x\omega) \rket^{\perp,m+1} \setminus \rbra T_{j^1_x\omega}S_\LL / \ker\widetilde\Omega(j^1_x\omega) \rket \subset T_{j^1_x\omega}J^1(\pi\circ\nu) / \ker\widetilde\Omega(j^1_x\omega)$ is not empty and $v$ is a vector belonging to this set, by linearity, we may assume that $v$ is a non-zero vector in $\abra\cbra\pp{p^i_\alpha}, \pp{p^i_{\alpha i}}\cket\aket$. Suppose that $v=v^i_\alpha\pp{p^i_\alpha}+v_\alpha\pp{p^i_{\alpha i}}$
\[ 0 = \tilde\Omega\Bigl(v,U^i_\alpha,X_1,\dots,X_m\Bigr) = -v^i_\alpha \quand
   0 = \tilde\Omega\Bigl(v,U_\alpha  ,X_1,\dots,X_m\Bigr) = -v_\alpha \,, \]
So $v$ is null, which is a contradiction.
\end{proof}

Next, we present some examples.

\begin{example}[\textbf{Affine Lagrangian densities}] \label{ex:lagrangian:affine}
Let $\gamma\colon E\to\MM\pi\cong(J^1\pi)^\dag$ be a section of the fiber bundle $\nu\colon\MM\pi\to E$ and consider the affine Lagrangian density given by its action on $J^1\pi$, that is,
\[ \LL_\gamma(j^1_x\phi) := \abra\gamma(\phi(x)),j^1_x\phi\aket\eta \,. \]
In adapted coordinates, if $\gamma=(x^i,u^\alpha,\gamma_0(x,u),\gamma^i_\alpha(x,u))$ and $j^1_x\phi=(x^i,u^\alpha,u^\alpha_i)$, then
\[ \LL_\gamma(x^i,u^\alpha,u^\alpha_i) = \rbra\gamma_0(x,u) + \gamma^i_\alpha(x,u)u^\alpha_i\rket\dmx \,. \]
Note that this Lagrangian is obviously degenerate since its Hessian velocity matrix is null.

Examples of this type of Lagrangians appear, for instance, in the metric-affine gravitation theory with the Hilbert-Einstein Lagrangian density. The Lagrangian density of Dirac fermion fields in the presence of a background tetrad field and a background spin connection is also affine (for more details, see \cite{GchtMngtSrdn97}).

The $(m+1)$--Lagrangian submanifold  $S_{\LL_\gamma} = (A_\pi^\eta)^{-1}(\diff\LL_\gamma(J^1\pi))$ of the premultisymplectic manifold $(J^1(\pi \circ \nu), \widetilde\Omega)$ is then given from the local expression \eqref{Apieta} of $A_\pi^\eta$ by
\[ S_{\LL_\gamma}= \cbra (x^i, u^\alpha, p, p^i_\alpha, u^\alpha_j, p_j, p^i_{\alpha j}) \,\colon\, p^i_\alpha=\gamma^i_\alpha \, ,\ p^i_{\alpha i}=\pp[\gamma_\circ]{u^\alpha}+\pp[\gamma_\beta^j]{u^\alpha}u^\beta_j \cket \,. \]
\hfill$\lhd$
\end{example}

\begin{example}[\textbf{Quadratic Lagrangian densities}] \label{ex:lagrangian:quadratic}
Let $\flat\colon J^1\pi\to\MM\pi\cong(J^1\pi)^\dag$ be a morphism of affine bundles over the identity of $E$. Then we define the quadratic Lagrangian function
\[ L_\flat(z) := \frac12\abra\flat(z),(z)\aket \,,\quad \forall z\in J^1\pi \,. \]
If we suppose that $\flat$ is locally given by
\[\flat(x^i, u^\alpha, u^\alpha_i) = (x^i, u^\alpha, \  \flat_\circ(x,u) + \flat^i_\alpha(x,u)u^\alpha_i, \ \tilde\flat^i_\alpha(x,u) + \flat^{ij}_{\alpha \beta}(x,u)u^\beta_j),\]
then the Lagrangian function $L_\flat$  is locally given by
\[ L_\flat(x^i, u^\alpha, u^\alpha_i) = \frac12\rbra \flat_\circ + (\flat^i_\alpha + \tilde{\flat}^i_\alpha)u^\alpha_i + \flat^{ij}_{\alpha\beta}u^\alpha_iu^\beta_j \rket \,. \]
Note that, in this case, the Hessian velocity matrix of $L_\flat$ is $\frac12(\flat^{ij}_{\alpha\beta}+\flat^{ji}_{\beta\alpha})$. Therefore, if $\flat$ is assumed to be symmetric ($\flat^{ij}_{\alpha\beta}=\flat^{ji}_{\beta\alpha}$), then $L_\flat$ will be regular if and only if the composition $\mu \circ \flat: J^1\pi \to J^1\pi^\circ$ is an affine bundle isomorphism.

Examples of this type of Lagrangian functions appear, for instance, in the theory of electromagnetic fields and Proca fields (refer again to \cite{GchtMngtSrdn97}).

The $(m+1)$--Lagrangian submanifold  $S_{\LL_\flat} = (A_\pi^\eta)^{-1}(\diff\LL_\flat(J^1\pi))$  of the premultisymplectic manifold $(J^1(\pi \circ \nu), \widetilde\Omega)$, where $\LL_\flat=L_\flat\eta$, is then given from the local expression \eqref{Apieta} of $A_\pi^\eta$ by
\begin{multline*}
S_{\LL_\flat}= \Bigg\{ (x^i, u^\alpha, p, p^i_\alpha, u^\alpha_j, p_j, p^i_{\alpha j}) \,\colon \,  p^i_\alpha=\frac12\rbra\flat^{ij}_{\alpha\beta}+\flat^{ji}_{\beta\alpha}\rket u^\beta_j \,, \\
p^i_{\alpha i}= \frac12\rbra\pp[\flat_\circ]{u^\alpha}+\pp[(\flat^i_\beta + \tilde\flat^i_\beta)]{u^\alpha}u^\beta_i +\pp[\flat^{ij}_{\beta \gamma}]{u^\alpha}u^\beta_iu^\gamma_j\rket \Bigg\} \,.
\end{multline*}
\hfill$\lhd$
\end{example}

We return to the general case.

\begin{theorem}
Given a Lagrangian density $\LL$, we have that:
\begin{enumerate}
\item A (local) section $\sigma\in\Sections(\pi)$ is a solution of the Euler-Lagrange equations if and only if
\[(A_\pi^\eta)^{-1}\circ \diff\LL \circ j^1\sigma = j^1(\Leg_\LL \circ j^1\sigma).\]
\item The local equations defining $S_\LL$ as an ($m+1$)--Lagrangian submanifold of $J^1(\pi\circ\nu)$ are just the Euler-Lagrange equations for $\LL$.
\end{enumerate}
\end{theorem}

\begin{proof}
A local computation, using \eqref{eq:euler-lagrange:coord}, \eqref{eq:legendre:transform:ext:coord} and \eqref{Apieta}, proves the result.
\end{proof}

\begin{figure}[h]
\[ \xymatrix{
&S_\LL\ar[dr]& \\
\Lambda^{m+1}_2(J^1\pi)\ar[rd]&&J^1(\pi \circ \nu)\ar[ll]_{A_\pi^\eta}\ar[ld]\ar[rd]&& \\
& J^1\pi\ar@<1ex>[ul]^{\diff\LL}\ar[rr]^{\Leg_\LL}\ar[rd]^{\pi_{1,0}}\ar[rdd]^{\pi_1}&&\MM\pi \cong J^1\pi^\dag\ar[ld]
_\nu\ar[ldd]_{\pi \circ \nu}&& \\
&&E\ar[d]^\pi&&& \\
&&M \ar `^u[rr] `[uuu]_{j^1(\Leg_\LL\circ j^1\sigma)} [uuu]
    \ar[u]<1ex>^\sigma\ar@/^1pc/[uul]^{j^1\sigma}\ar@/_1pc/[uur]_{\Leg_\LL\circ j^1\sigma}&&&
} \]
\caption{The Lagrangian formalism in the Tulczyjew's triple}
\label{fig:t3cft:lagrangian}
\end{figure}
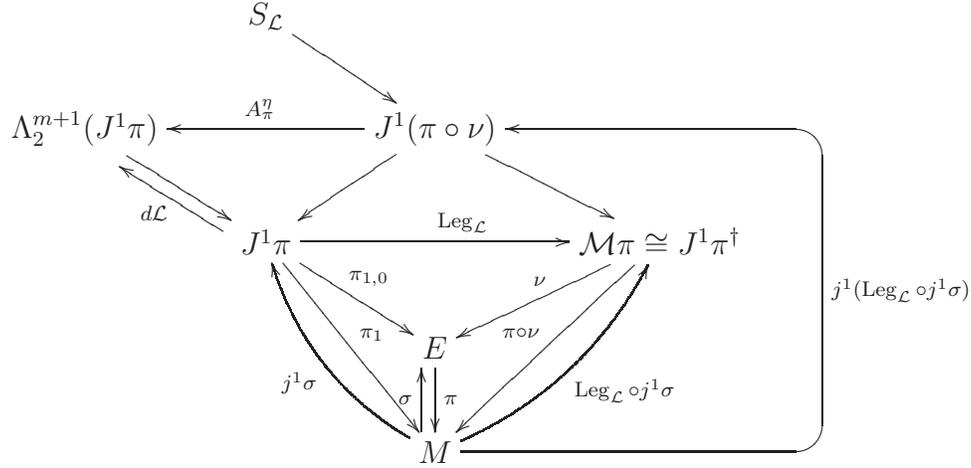

Figure \ref{fig:t3cft:lagrangian} illustrates the above situation

%%% HAMILTONIAN FORMALISM ------------------------------------------------------
\subsection{Hamiltonian formalism}
The construction of the Hamiltonian part of the Tulczyjew's triple for Classical Field Theory involves the 1st-jet bundle $J^1(\pi\circ\nu)$, which has already been introduced in the previous section, and the bundle of forms $\Lambda^{m+1}_2\MM\pi$. An arbitrary element $\bar\omega$ of $\Lambda^{m+1}_2\MM\pi$ is locally written
\[ \bar p_\alpha\du^\alpha\wedge \dmx + \bar p\diff p\wedge \dmx + \bar p_i^\alpha \diff p^i_\alpha\wedge \dmx \,, \]
so we may consider natural local coordinates $(x^i, u^\alpha, p, p^i_\alpha, \bar p_\alpha, \bar p, \bar p^\alpha_i)$ on $\Lambda^{m+1}_2\MM\pi$. We recall that, as usual, adapted coordinates on $J^1(\pi\circ\nu)$ are denoted $(x^i, u^\alpha, p, p_{ \alpha}^i, u^\alpha_j, p_j, p^i_{ \alpha j})$.

In general, for any fiber bundle $\pi\colon E\to M$, a jet $j^1_x\phi\in J^1\pi$ defines a horizontal projector in $T_{\phi(x)}E$, $\hor_{j^1_x\phi}:=T_x\phi\circ T_{\phi(x)}\pi$. In local coordinates,
\[ \hor_{j^1_x\phi} = \dx^i\tensor\rbra \pp{x^i} + u^\alpha_i\pp{u^\alpha} \rket \,. \]
In the particular case of $\pi\circ\nu\colon\MM\pi\to M$,
\[ \hor_{\bar z} = \dx^j\tensor\rbra \pp{x^j} + u^\alpha_j\pp{u^\alpha} + p_j\pp{p} + p^i_{\alpha j}\pp{p^i_\alpha} \rket \,, \]
for any $\bar z=(x^i, u^\alpha, p, p_{ \alpha}^i, u^\alpha_j, p_j, p^i_{ \alpha j})\in J^1(\pi\circ\nu)$.

We define the affine bundle morphism $\flat_\Omega\colon J^1(\pi\circ\nu) \to \Lambda^{m+1}_2(\MM\pi)$
\begin{equation} \label{empb}
\flat_\Omega(\bar z) = i_{\hor_{\bar z}}\Omega(\omega) - (m-1)\Omega(\omega)
\end{equation}
where $\omega=(\pi\circ\nu)_{1,0}(\bar z)$, $\Omega$ is the canonical multisymplectic form of $\MM\pi$ and
\[\rbra i_{\hor_{\bar z}}\Omega(\omega)\rket (X_1,\ldots, X_{m+1}) = \sum_{i=1}^{m+1}\Omega(\omega) (X_1,\ldots, \hor_{\bar z}(X_i),\ldots, X_{m+1})\]
for $X_1,\ldots,X_{m+1} \in T_\omega\MM\pi$. Locally, we have that
\begin{equation} \label{eq:flat:omega:coord}
\flat_\Omega(x^i, u^\alpha, p, p_{ \alpha}^i, u^\alpha_j, p_j, p^i_{ \alpha j}) = \rbra x^i, u^\alpha, p, p_{ \alpha}^i, p_{\alpha j}^j, - 1 , - u^\alpha_j\rket \,.
\end{equation}

According to Example \ref{ex:multisymplectic:horizontal}, we know that $\Lambda^{m+1}_2\MM\pi$ has a canonical multisymplectic structure. Namely,
\begin{equation} \label{eq:Omega:LambdaMpi:coord}
\Omega_{\Lambda^{m+1}_2\MM\pi} = -\diff\bar p_\alpha\wedge\du^\alpha\wedge \dmx - \diff\bar p\wedge \diff p\wedge \dmx - \diff\bar p_i^\alpha\wedge \diff p^i_\alpha\wedge \dmx \,.
\end{equation}
We pullback this structure by $\flat_\Omega$, obtaining a premultisymplectic ($m+2$)--form on $J^1(\pi\circ\nu)$. From Equations \eqref{eq:flat:omega:coord} and \eqref{eq:Omega:LambdaMpi:coord}, it follows that
\begin{equation} \label{eq:omega:tilde:coord:alt}
\flat_\Omega^*(\Omega_{\Lambda^{m+1}_2\MM\pi}) = -\diff p^j_{\alpha j}\wedge\du^\alpha\wedge \dmx - \diff p^i_\alpha\wedge\du^\alpha_i\wedge \dmx \,.
\end{equation}
Thus, $\flat^*_\Omega(\Omega_{\Lambda^{m+1}_2\MM\pi})$ turns out to be the ($m+2$)--form $\widetilde\Omega$ on $J^1(\pi\circ\nu)$ introduced in Equation \eqref{eq:omega:tilde:coord} of Section \secref{sec:t3cft:lagrangian}.

\begin{theorem} \label{th:tulczyjew:morphism:right}
The morphism $\flat_\Omega\colon J^1(\pi\circ\nu) \to \Lambda^{m+1}_2(\MM\pi)$, locally given by Equation \eqref{eq:flat:omega:coord}, is an affine bundle epimorphism over the identity of $\MM\pi$. Moreover, it is canonical in the sense that it only depends on the original fiber bundle $\pi\colon E\to M$.
\end{theorem}

Let $\HH\colon\MM\pi\to\Lambda^mM$ be a Hamiltonian density, then $\diff\HH\colon\MM\pi\to\Lambda^{m+1}_2E$. Using \eqref{eq:hamiltonian:density:coord} and \eqref{eq:flat:omega:coord}, it follows that
\[ -\diff\HH(\MM\pi) \subseteq \flat_\Omega(J^1(\pi\circ\nu)) \,. \]
In fact, if $\HH(x^i, u^\alpha, p, p^i_\alpha) = (p+H(x^i, u^\alpha, p^i_\alpha))\dmx$, we deduce that
\begin{equation} \label{eq:hamiltonian:differential}
-\diff\HH = - \pp[H]{u^\alpha}\du^\alpha\wedge \dmx -  \diff p\wedge \dmx- \pp[H]{p_\alpha^j}\diff p_\alpha^j\wedge \dmx \,.
\end{equation}

\begin{proposition}
$S_\HH = \flat_\Omega^{-1}\rbra-\diff\HH(\MM\pi)\rket$ is an ($m+1$)--Lagrangian submanifold of the premultisymplectic manifold $(J^1(\pi\circ\nu), \widetilde\Omega)$.
\end{proposition}

\begin{proof}
The proof is analogous to that of Proposition \ref{th:SL:lagrangian:manifold}. Nonetheless, it is worth noting that, from \eqref{eq:flat:omega:coord} and \eqref{eq:hamiltonian:differential}, we obtain that the submanifold $S_\HH$ is locally given by

\begin{equation} \label{eq:Sh}
S_\HH = \cbra(x^i, u^\alpha, p, p^i_\alpha, u^\alpha_j, p_j, p^i_{\alpha j}):  u^\alpha_j = \pp[H]{p_\alpha^j} \,, \quad p^j_{\alpha j} = -\pp[H]{u^\alpha} \cket \,.
\end{equation}
\end{proof}

Next, we present an example.

\begin{example}[\textbf{Quadratic Hamiltonian densities}]
Let $\sharp\colon J^1\pi^\circ\cong\MM^\circ\pi \to J^1\pi$ be a morphism of affine bundles over the identity of $E$. Then, we may define the Hamiltonian density $\HH_\sharp\colon \MM\pi\cong J^1\pi^\dag \to \Lambda^mM$ given by
\[ \HH_\sharp(\omega) := \abra\omega,\sharp(\mu(\omega))\aket\eta \,,\quad \forall\omega\MM\pi \,. \]
If $\sharp$  is locally given by
\[ \sharp(x^i, u^\alpha, p^i_ \alpha) = (x^i, u^\alpha, \sharp^\alpha_i(x,u) + \sharp^{\alpha\beta}_{ij}(x,u) p^j_\beta )\]
then $\HH_\sharp$ is locally given by
\[ \HH_\sharp(x^i, u^\alpha, p, p^i_ \alpha) =  \rbra p + \sharp^\alpha_i(x,u)p^i_\alpha + \sharp^{\alpha\beta}_{ij}(x,u) p^i_\alpha p^j_\beta \rket\dmx \,,\]
with associated Hamiltonian function
\[ H_\sharp(x^i, u^\alpha, p^i_ \alpha) =  \sharp^\alpha_i(x,u)p^i_\alpha + \sharp^{\alpha\beta}_{ij}(x,u) p^i_\alpha p^j_\beta \,.\]
From these expressions, we see that the Hamiltonian $\HH_\sharp$ is quadratic.

Examples of this type of Hamiltonian sections may be obtained from hyper-regular quadratic Lagrangian functions (see Example \ref{ex:lagrangian:quadratic}).

The $(m+1)$--Lagrangian submanifold $S_{\HH_\sharp} = \flat_\Omega^{-1} \rbra -\diff\HH_\sharp(\MM\pi)\rket$ of the premultisymplectic manifold $(J^1(\pi \circ \nu), \widetilde\Omega)$ is then given from the local expression \eqref{eq:flat:omega:coord} of $\flat_\Omega$ by
\begin{multline*}
S_{\HH_\sharp}= \Bigg\{ (x^i, u^\alpha, p, p^i_\alpha, u^\alpha_j, p_j, p^i_{\alpha j}) \,\colon\, u^\alpha_ i=  \sharp^\alpha_i(x,u) + \rbra \sharp_{ij}^{\alpha\beta}(x,u)+ \sharp_{ji}^{\beta\alpha}(x,u) \rket p_\beta^j \,, \\
\ p^i_{\alpha i}= - \pp[\sharp_j^\beta(x,u)]{u^\alpha}p_\beta^j - \pp[\sharp^{\beta\gamma}_{jk}(x,u)]{u^\alpha}p^j_\beta p^k_\gamma \Bigg\} \,.
\end{multline*}
\hfill$\lhd$
\end{example}

We return to the general case.
\begin{theorem}
Given a Hamiltonian section $h\in\Sections(\mu)$, let $\HH\colon \MM\pi \to \Lambda^mM$ be the associated Hamiltonian density. We have that
\begin{enumerate}
\item A section $\tau\colon M\to \MM^\circ\pi$ is a solution of the Hamilton-De Donder-Weyl equation if and only if
\[\flat_\Omega \circ j^1(h\circ\tau) = -\diff\HH \circ (h\circ\tau).\]
\item The local equations defining $S_\HH$ as an ($m+1$)--Lagrangian submanifold of $J^1(\pi\circ\nu)$ are just the Hamilton's equations.
\end{enumerate}
\end{theorem}

\begin{proof}
A local computation, using \eqref{eq:hamiltonian:density:coord} and \eqref{eq:flat:omega:coord}, proves the result.
\end{proof}

\begin{figure}[h]
\[
\xymatrix{
&S_\HH\ar[dl]& \\
J^1(\pi\circ\nu)\ar[dr]_{(\pi\circ\nu)_{1,0}}\ar[rr]^{\flat_\Omega}&&\Lambda^{m+1}_2\MM\pi\ar[dl] \\
&\MM\pi\ar[dddl]_{\pi\circ\nu}\ar@<1ex>[d]^\mu\ar@<-1ex>[ur]_{-\diff\HH}\ar[ddr]^\nu& \\
&\MM^\circ\pi\ar[u]^h\ar[dr]_{\pi^*_{1,0}}& \\
&&E\ar[dll]^\pi \\
M\ar[uur]_\tau\ar[uuuu]^{j^1(h\circ \tau)}&&
}
\]
\caption{The Hamiltonian formalism in the Tulczyjew's triple}\label{s}
\end{figure}
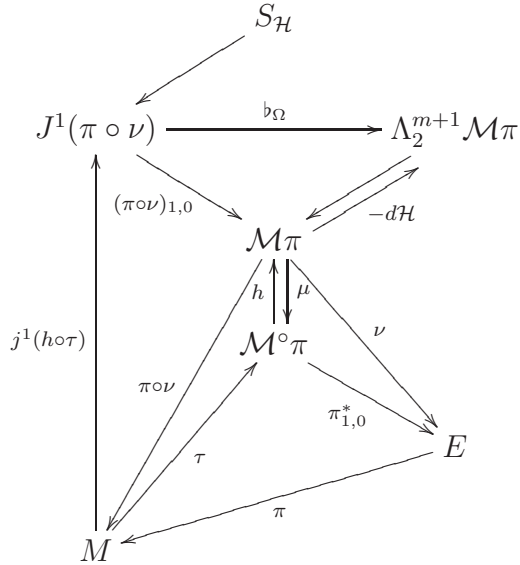

Figure \ref{s} illustrates the above situation

%%% EQUIVALENCE BETWEEN BOTH FORMALISMS ----------------------------------------
\subsection{Equivalence between both formalisms}
In section \secref{sec:cft:equivalence}, we already saw how the Lagrangian and Hamiltonian formalisms are interrelated by means of the Legendre transform. Here, we recover this relation within the Tulczyjew's triple for Classical field theory by pulling the dynamical structures of $\Lambda^{m+1}_2J^1\pi$ and $\Lambda^{m+1}_2\MM\pi$ to $J^1(\pi\circ\nu)$ as Lagrangian submanifolds.

Let $\LL\colon J^1\pi\to\Lambda^mM$ be an hyperregular Lagrangian density and $\HH$ the associated Hamiltonian density (see Equation \eqref{eq:legendre:hamiltonian:density}). A simple computation in local coordinates using Equations \eqref{eq:SL}, \eqref{eq:Sh} and \eqref{eq:legendre:hamiltonian:function:coord} shows the following result.

\begin{theorem}
The ($m+1$)--Lagrangian submanifolds $S_\LL = (A_\pi^\eta)^{-1}(\diff\LL(J^1\pi))$  and $S_\HH= \flat_\Omega^{-1}(-\diff\HH(\MM\pi))$ of the premultisymplectic manifold $(J^1(\pi\circ\nu), \widetilde\Omega)$ are equal.
\end{theorem}

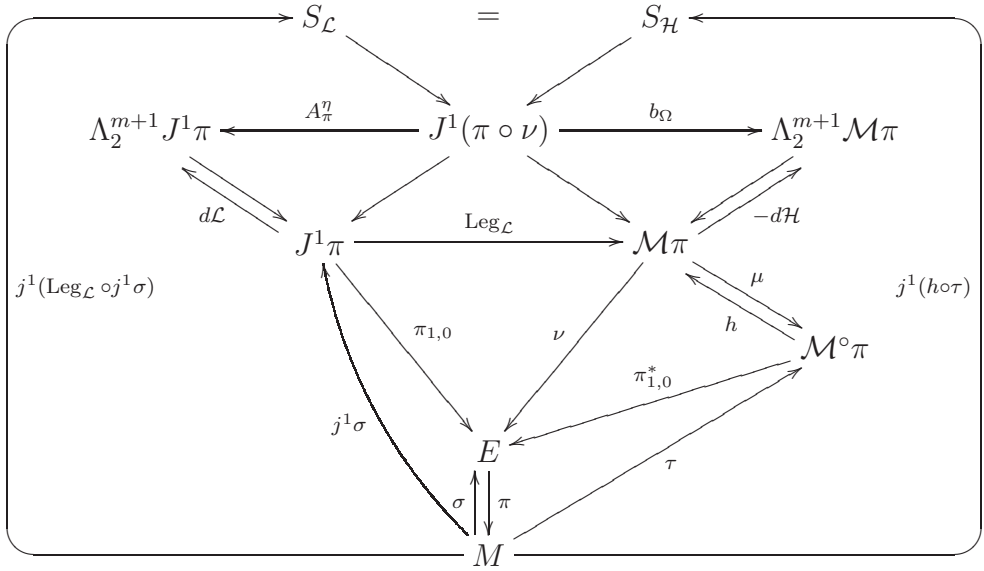
\begin{figure}[h]
\[
\xymatrix{
&&S_\LL\ar[dr]&=&S_\HH\ar[dl]&&& \\
&\Lambda^{m+1}_2J^1\pi\ar[rd]&&J^1(\pi\circ\nu)\ar[rr]^{b_\Omega}\ar[ll]_{A_\pi^\eta}\ar[ld]\ar[rd]&&\Lambda^{m+1}_2\MM\pi\ar[dl]& \\
&& J^1\pi\ar@<1ex>[ul]^{\diff\LL}\ar[rr]^{\Leg_\LL}\ar[rdd]^{\pi_{1,0}}&& \MM\pi\ar@<-1ex>[ur]_{-\diff\HH }\ar[dr]^\mu\ar[ldd]_\nu&&& \\
&&&&&\MM^\circ\pi\ar@<1ex>[ul]^h\ar[dll]_{\pi_{1,0}^*}& \\
&&&E\ar[d]^\pi&&&& \\
&&&M\ar[uurr]_\tau\ar[u]<1ex>^\sigma\ar@/^1pc/[uuul]^{j^1\sigma}
\ar `^u[rrr] `[uuuuu]^{j^1(h\circ \tau)} [ruuuuu]
\ar `_u[lll] `[uuuuu]_{j^1(\Leg_\LL\circ j^1\sigma)} [luuuuu]&&&
}
\]

\caption{The Tulczyjew's triple for Classical Field Theories}\label{k}
\end{figure}

Figure \ref{k} illustrates this situation.

% CONCLUSIONS AND FUTURE WORK --------------------------------------------------
\section{Conclusions and future work}
A Tulczyjew's triple associated with a fibration is introduced. This construction allows us to describe Euler-Lagrange and Hamilton-De Donder-Weyl equations for classical field theories as the local equations defining Lagrangian submanifolds of a premultisymplectic manifold.

It would be interesting to extend these results for the more general case of classical field theory on Lie algebroids using a multisymplectic formalism. The first steps have been done in this direction (see \cite{GzmMrrVnk12}). We remark that these ideas could be applied to the reduction of classical field theories  which are invariant under the action of a symmetry Lie group.

Another interesting goal is to develop a geometric formulation of vakonomic (nonholonomic) classical field theories using our Tulczyjew's triple. Extensions of this construction for vakonomic (nonholonomic) classical theories with symmetries (using the Lie algebroid setting) could also be discussed.

% APPENDIX ---------------------------------------------------------------------
\appendix

% VOLUME INDEPENDENT FORMULATION OF T3 -----------------------------------------
\section{Volume independent formulation of the Tulczyjew's triple}
\label{sec:t3cft:novol}
Even though the Euler-Lagrange (and Hamilton-De Donder-Weyl) equations of Classical Field Theory are usually given using a volume form on the base space of the configuration bundle, they can be obtained without the need of that volume form. In fact, one may establish such equations in the more general case where the base manifold is not necessarily orientable (see, for instance, \cite{Cmps10,CmpsLnMrtVnk09,GchtMngtSrdn97}). %A volume form ensures the existence of globally defined non-degenerate solutions, but there are physical examples where there aren't.

In this sense, it is important to be able to give a Tulczyjew's triple independent of any chosen basic volume form. This is the aim of this Appendix. To do so, we construct the left Tulczyjew's map in three steps: First, we transform the 1st-jet bundle of $\Lambda^m\pi_M\colon\Lambda^mM\to M$ into a ``simpler'' space. Second, we compute (in some sense) the core of the Tulczyjew's morphism from the ``tangent'' pairing and the involution. Third, we combine the previous two steps in order to obtain the desired map. Before we start, we present a short review of the multisymplectic formulation of Classical Field Theory without the need of the orientability assumption (compare with Section \secref{sec:cft:multisymplectic}). Along the rest of the paper, local coordinates on $M$ will no longer be required compatible with any volume form.

\ \\\paragraph{\textbf{Volume independent formulation of Classical Field Theory:}} Given a Lagrangian density $\LL\colon J^1\pi \to \Lambda^mM$, the associated Lagrangian function $L$ is now locally defined such that $\LL=L\dmx$. Of course, in presence of a volume form $\eta$ on $M$ (if such is the case), it will coincide on a compatible chart with the globally defined Lagrangian function $\bar L\colon J^1\pi \to \RR$ such that $\LL=\bar L\eta$.

In this setting, the Poincaré-Cartan $m$--form is given by the expression
\begin{equation} \label{eq:pc-mform}
\Theta_\LL = \LL + \abra S,\diff\LL \aket \,,
\end{equation}
where $S\colon TJ^1\pi \to \pi_1^*(TM) \tensor_{J^1\pi} \Vert(\pi_{1,0})$ is the canonical vertical endomorphism, which has the local representation
\begin{equation} \label{eq:vertical:endomorphism:canonical}
S = \rbra\du^\alpha-u^\alpha_j\dx^j\rket \tensor \pp{x^i} \tensor \pp{u^\alpha_i} \,.
\end{equation}
This newly defined Poincaré-Cartan $m$--form coincides in fact with the previously defined one, Equation \eqref{eq:pc-forms}, and therefore the Euler-Lagrange equations \eqref{eq:euler-lagrange:coord} follow.

In the other side of the picture, we redefine the dual bundles of $J^1\pi$. Instead of considering as extended dual of $J^1\pi$ the affine maps on $J^1\pi$ with values in $\RR$, we consider the affine maps with values in $\Lambda^mM$. The reduced dual of $J^1\pi$ will be the quotient of the extended dual by fiberwise constant affine maps on $J^1\pi$ with values in $\Lambda^mM$. Resuming,
\begin{align}
 (J^1\pi)^\dag &:= \Aff_\pi(J^1\pi,\Lambda^mM)  = \cbra A\in\Aff(J^1_u\pi,\Lambda^m_{\pi(u)}M) \,\colon\, u\in E \cket \,, \\
(J^1\pi)^\circ &:= \Aff_\pi(J^1\pi,\Lambda^mM) / \mathrm{Morph}_M(E,\Lambda^mM) \,.
\end{align}

The improvement of considering $\Lambda^mM$ as the target space of the affine maps is that the identifications
\begin{equation} \label{eq:dual:identification:canonical}
J^1\pi^\dag \cong \MM\pi \quand J^1\pi^\circ \cong \MM^\circ\pi
\end{equation}
are now canonical, in contrast with those of \eqref{eq:dual:identification}, which depend on a volume form on the base. In order to emphasize this and avoid confusion with Section \secref{sec:t3cft}, we will use the dual affine bundle notation $J^1\pi^\dag$ and $J^1\pi^\circ$ in spite of the the form bundle one $\MM\pi$ and $\MM^\circ\pi$.

The overall notation and objects (coordinates, forms, etc.) will remain unaltered, but the natural pairing is now slightly different:
\[ \abra(x^i, u^\alpha, p, p_\alpha^i),\, (x^i, u^\alpha, u^\alpha_i)\aket = (p+p_\alpha^iu^\alpha_i)\dmx \,. \]

Another technical issue due to the lack of the volume form is that, the line bundle $\mu\colon J^1\pi^\dag \to J^1\pi^\circ$ cannot be seen as an $\RR$--principal bundle. Therefore, the definition of a Hamiltonian density $\HH\colon J^1\pi^\dag \to \Lambda^mM$ must be reformulated by the condition
\[ i_V(\Omega+\diff\HH) = 0 \,,\quad \forall V\in\Vert(\mu) \,. \]
From here, Hamilton's equations \eqref{eq:hamilton:coord} follow immediately.

\ \\\paragraph{\textbf{The 1st-jet bundle of the bundle of $m$--forms:}} The 1st-jet bundle $J^1\pi$ of an arbitrary fiber bundle $\pi\colon E\to M$ is, in general, only an affine bundle over $E$. But, in some particular cases, it can be can be seen as a vector bundle over $M$; for instance when $E$ is a tensor bundle of $M$, \ie\ $E=TM,T^*M,\Lambda^mM,\dots$, and $M$ itself is provided with a linear connection $\nabla$. For obvious reasons, we expose the particular case of the tensor bundle of $m$--forms $\Lambda^mM$.

If we think of the elements of $J^1\Lambda^m\pi_M$ as linear maps from $TM$ to $T(\Lambda^m\pi_M)$ (see Section \ref{sec:cft:lagrangian:formalism}), we may define the bundle isomorphism
\[ \begin{array}{rcl}
      J^1\Lambda^m\pi_M \subset (\Lambda^m\pi_M)^*(T^*M)\tensor_{\Lambda^mM}T(\Lambda^m\pi_M) &\To   & (\Lambda^m\pi_M)^*(T^*M)\tensor_{\Lambda^mM}\Vert(\Lambda^m\pi_M)\\
      j^1_x\alpha = T_x\alpha &\Mapsto& j^1_x\alpha - j^1\tilde\alpha = T_x\alpha - T_x\tilde\alpha \,,
\end{array} \]
where $\tilde\alpha$ is any (local) section of $\Lambda^m\pi_M\colon\Lambda^mM\to M$ such that
\[ \tilde\alpha(x)=\alpha(x) \quand (\nabla_X\tilde\alpha)(x)=0 \,, \forall X\in\Fields(M) \,. \]
It is easy to prove that the above defined morphism does not depend on the chosen section $\tilde\alpha$. In fact, if $\alpha(x) = a\dmx|_x$ $\tilde\alpha=\tilde a\dmx$, the previous conditions imply that $\pp*[\tilde a]{x^i}|_x=a\cdot\chr^j_{ij}(x)$, which determines the 1st-jet of $\tilde\alpha$. Therefore, for local coordinates $(x^i,a,a_i)$ on $J^1\Lambda^m\pi_M$, we have
\[  j^1_x\alpha-j^1_x\tilde\alpha = \rbra x^i,a,a_i-a\cdot\chr^j_{ij}\rket = \rbra a_i-a\cdot\chr^j_{ij}\rket\cdot\dx^i\tensor\pp{a} \,. \]

Besides of this, we have that the vertical bundle $\Vert(\Lambda^m\pi_M)$ is isomorphic to $\Lambda^mM\times_M\Lambda^mM$ by the vertical lift. Namely,
\[ \begin{array}{rcl}
    (\cdot)^\vert_\cdot\colon
    \Lambda^mM\times\Lambda^mM &\To    & \Vert(\Lambda^m\pi_M) \\
                  (\alpha,\beta) &\Mapsto& (\beta)^\vert_\alpha \\
       (x^i,a,b) = (a\dmx,b\dmx) &\Mapsto& (x^i,a,b) = b\pp{a}\big|_{a\dmx} \,.
\end{array} \]
Therefore, we have for the the vector bundle associated to $J^1\Lambda^m\pi_M$
\[ \begin{array}{rcl}
(\Lambda^m\pi_M)^*(T^*M)\tensor\Vert(\Lambda^m\pi_M)
&\cong& (\Lambda^m\pi_M)^*(T^*M)\tensor(\Lambda^mM\times\Lambda^mM) \\
&\cong& \Lambda^mM\times(T^*M\tensor\Lambda^mM) \,,
\end{array} \]
which in local coordinates reads
\[ \begin{array}{rcl}
(\Lambda^m\pi_M)^*(T^*M)\tensor\Vert(\Lambda^m\pi_M) &\To    & \Lambda^mM\times(T^*M\tensor_M\Lambda^mM) \\
                                                     (x^i,a,a_i) &\Mapsto& (x^i,a,a_i) \\
                        a_i\dx^i\tensor\pp{a}\big|_{a\dmx} &\Mapsto& (a\dmx, a_i\dx^i\tensor\dmx) \,.
\end{array} \]

Summing up, we obtain an affine bundle transformation $\Phi^\nabla$ over $\Lambda^mM$ from $J^1\Lambda^m\pi_M$ to $\Lambda^mM\times_M(T^*M\tensor\Lambda^mM)$ which depends on the linear connection $\nabla$. Locally, we have
\[ \begin{array}{rcl}
         \Phi^\nabla \colon J^1\Lambda^m\pi_M &\To    & \Lambda^mM\times_M(T^*M\tensor\Lambda^mM) \\
                                  (x^i,a,a_i) &\Mapsto& (x^i,a,a_i-a\cdot\chr^j_{ij}) \\
\dx^i\tensor(\pp{x^i}+a_i\pp{a}\big|_{a\dmx}) &\Mapsto& (a\dmx, (a_i-a\cdot\chr^j_{ij})\dx^i\tensor\dmx) \,.
\end{array} \]

\ \\\paragraph{\textbf{The core map:}} The natural pairing between the elements of $J^1\pi^\dag$ and $J^1\pi$ is the fibered map
\[ \begin{array}{rcl}
               \pairing\colon J^1\pi^\dag\times_EJ^1\pi &\To    & \Lambda^mM \\
                                     (\omega,j^1_x\phi) &\Mapsto& \abra\omega,j^1_x\phi\aket \\
((x^i,u^\alpha,p,p^i_\alpha),(x^i,u^\alpha,u^\alpha_i)) &\Mapsto& (p+p^i_\alpha u^\alpha_i)\dmx \,.
\end{array} \]
We lift it by 1st-prolongation to the pairing
\[ \begin{array}{rcl}
j^1\pairing\colon
J^1\pi^\dag\times_{J^1\pi}(J^1\pi_1,j^1(\pi_{1,0}),J^1\pi) &\To    & J^1\Lambda^m\pi_M \\
                                 (j^1_x\omega,j^1_x\sigma) &\Mapsto& j^1(\abra\omega,\sigma\aket) \\
((x^i,u^\alpha,p,p^i_\alpha,\bar u^\alpha_j,p_j,p^i_{\alpha j}),(x^i,u^\alpha,u^\alpha_i,\bar u^\alpha_j,u^\alpha_{ij})) &\Mapsto& (x^i,p+p^i_\alpha u^\alpha_i,p_k+p^i_{\alpha k}u^\alpha_i+p^i_\alpha u^\alpha_{ik}) \,,
\end{array} \]
where the second argument fibers over $J^1\pi$ by $j^1(\pi_{1,0})$. Composing it with the exchange transformation, we get a map whose coordinate representation is
\[ \begin{array}{rcl}
j^1\abra\cdot,\ex_\nabla(\cdot)\aket\colon
J^1(\pi^\dag_1)\times_{J^1\pi}(J^1\pi_1,(\pi_1)_{1,0},J^1\pi) &\To    & J^1\Lambda^m\pi_M \\
                                                   (j^1_x\omega,j^1_x\sigma) &\Mapsto& j^1\abra j^1_x\omega,\ex_\nabla(j^1_x\sigma)\aket \\
((x^i,u^\alpha,p,p^i_\alpha,u^\alpha_i,p_j,p^i_{\alpha j}),(x^i,u^\alpha,u^\alpha_i,\bar u^\alpha_j,u^\alpha_{ij})) &\Mapsto& (x^i,a,a_k) \,,
\end{array} \]
where now the second argument fibers over $J^1\pi$ by $(\pi_1)_{1,0}$, and where $a=p+p^i_\alpha \bar u^\alpha_i$ and $a_k=p_k+p^i_{\alpha k}\bar u^\alpha_i+p^i_\alpha(u^\alpha_{ki}+(\bar u^\alpha_j-u^\alpha_j)\chr^j_{ki})$.

\ \\\paragraph{\textbf{The Tulczyjew's map:}} Composing the previous transformations $\Phi^\nabla$ and $j^1\abra\cdot,\ex_\nabla(\cdot)\aket$ together with the 2nd component projection $\pr_2\colon\Lambda^mM\times_M(T^*M\tensor\Lambda^mM)\to T^*M\tensor\Lambda^mM$, we obtain the map
\[  \begin{array}{rcl}
\pr_2\circ\Phi^\nabla\circ j^1\abra\cdot,\ex_\nabla(\cdot)\aket\colon
J^1(\pi^\dag_1)\times_{J^1\pi}(J^1\pi_1,(\pi_1)_{1,0},J^1\pi) &\To    & T^*M\tensor\Lambda^mM \\
                                    (j^1_x\omega,j^1_x\sigma) &\Mapsto& j^1\abra j^1_x\omega,\ex_\nabla(j^1_x\sigma)\aket \\
((x^i,u^\alpha,p,p^i_\alpha,u^\alpha_i,p_j,p^i_{\alpha j}),(x^i,u^\alpha,u^\alpha_i,\bar u^\alpha_j,u^\alpha_{ij})) &\Mapsto& a_k\dx^k\tensor\dmx \,,
\end{array} \]
where $a_k=p_k+p^i_{\alpha k}\bar u^\alpha_i+p^i_\alpha(u^\alpha_{ki}+(\bar u^\alpha_j-u^\alpha_j)\chr^j_{ki}) - (p+p^i_\alpha \bar u^\alpha_i)\chr^j_{kj}$. Note that this map is a vector valued morphism that takes values from the affine bundles $(J^1(\pi^\dag_1),j^1(\pi^\dag_{1,0}),J^1\pi$ and $(J^1\pi_1,(\pi_1)_{1,0},J^1\pi)$. Hence it induces a vector valued affine bundle morphism over the identity of $J^1\pi$
\[ \widetilde{A_\pi^\nabla}\colon J^1(\pi^\dag_1) \To \Aff_{\pi_1}(J^1\pi_1,T^*M\tensor\Lambda^mM) \,, \]
where
\[ \Aff_{\pi_1}(J^1\pi_1, T^*M\tensor\Lambda^mM) = \cup_{z\in J^1\pi} \, \Aff(J^1_z\pi_1, T^*_{\pi_1(z)}M\tensor\Lambda^m_{\pi_1(z)}M) \,. \]
For adapted coordinates $(x^i, u^\alpha, u^\alpha_i, \bar p_k, \bar p^j_{\alpha k}, \bar p^{ij}_{\alpha k})$ on the vector bundle $\Aff_{\pi_1}(J^1\pi_1, T^*M\tensor\Lambda^mM)$,
\begin{multline*}
\widetilde{A_\pi^\nabla}(x^i, u^\alpha, p, p_\alpha^i, u^\alpha_j, p_j, p_{\alpha j }^i) =\\
= (x^i, u^\alpha, u^\alpha_i,
\bar p_k = p_k - p^i_\alpha u^\alpha_l\chr^l_{ki} - p\chr^j_{kj},
\bar p^i_{\alpha k} = p^i_{\alpha k} + p^l_\alpha\chr^i_{kl} - p^i_\alpha\chr^j_{kj},
\bar p^{ij}_{\alpha k} = p^j_\alpha\delta^i_k) \,.
\end{multline*}
By Lemma \ref{th:affine.map.lemma}, we have
\[ \Aff(J^1\pi_1,T^*M\tensor\Lambda^mM) \cong \Aff(J^1\pi_1,\Lambda^mM)\tensor T^*M \cong \Lambda^m_2J^1\pi\tensor T^*M \hookrightarrow \Lambda^{m+1}_2J^1\pi \,. \]
where the last ``inclusion'' is the natural morphism given by the wedge product, namely
\[ i\colon\omega\tensor\alpha\in\Lambda^m_2J^1\pi\tensor T^*M\Mapsto\alpha\wedge\omega\in\Lambda^{m+1}_2J^1\pi \,. \]
Therefore, if we assume $\nabla$ symmetric, we end up with the map
\[ \begin{array}{rcl}
A_\pi\colon J^1(\pi^\dag_1) &\To    & \Lambda^{m+1}_2J^1\pi \\
j^1_x\omega &\Mapsto& (i\circ\pr_2\circ\Phi^\nabla)(j^1\abra j^1_x\omega,\ex_\nabla(\cdot)\aket) \\
(x^i,u^\alpha,p,p^i_\alpha,u^\alpha_i,p_j,p^i_{\alpha j}) &\Mapsto& (p^i_{\alpha i}\du^\alpha+p^i_\alpha\du^\alpha_i)\wedge\dmx \,,
\end{array} \]
which we call \emph{the (left) Tulczyjew's morphism}.  It is important to note that, even though a symmetric linear connection $\nabla$ was necessary for the construction of $A_\pi$, actually the latter does not depend on the former. Hence, the Tulczyjew's morphism $A_\pi$ is a canonical transformation of $J^1(\pi^\dag_1)$.

% ON THE EXCHANGE MAP ----------------------------------------------------------
\section{On the exchange map $\ex_\nabla$}
The exchange map $\ex_\nabla\colon J^1\pi_1\to J^1\pi_1$ associated to a symmetric linear connection $\nabla$ on the base manifold $M$ of the fiber bundle $\pi\colon E\to M$ was first introduced by M. Modugno in \cite{Mdgn89} and, later on, developed in more detail together with his collaborator I. Kol\'a\v{r} in \cite{KlrMdgn91}. However, there is a small error in the coordinate expression of $\ex_\nabla$ in \cite{Mdgn89}, which also appears in \cite{KlrMdgn91}. Still, the notation in the later may lead to confusion and an unaware reader may fall in the same mistake. Since the exchange map is of crucial importance for the construction of the Tulczyjew's triple, we believe that is worth to clarify this matter here.

In \cite{Mdgn89}, Theorem 1 (page 359) establishes the coordinate expression of the exchange map. Following Modugno's notation:
\[ (x^\lambda,y^i,y^i_\mu,y^i_{\lambda0},y^i_{\lambda\mu})\circ s_k =
   (x^\lambda,y^i,y^i_{\lambda0},y^i_\mu,y^i_{\mu\lambda}+k^\nu_{\mu\lambda}(y^i_\nu-y^i_{\nu0})) \,, \]
where $s_k$ is the exchange maps associated to the linear connection $k$ with Christoffel symbols $k^\nu_{\lambda\mu}$. The connection term has a wrong sign and $s_k$ should be written
\[ (x^\lambda,y^i,y^i_\mu,y^i_{\lambda0},y^i_{\lambda\mu})\circ s_k =
   (x^\lambda,y^i,y^i_{\lambda0},y^i_\mu,y^i_{\mu\lambda}+k^\nu_{\mu\lambda}(y^i_{\nu0}-y^i_\nu)) \,. \]
In fact, the proof using a change of coordinates is incorrect and fails at the beginning of the page 360.

Another worth noting point is that, in order to define the exchange map $\ex_\nabla$, the connection $\nabla$ does not need to be symmetric. Although, in such a case, $\ex_\nabla$ is no longer an involution but an isomorphism and the construction of the Tulczyjew's triple must be carried out with a slight change. Instead of using the morphisms $\ex_\nabla$ and $\Phi^\nabla$ (or $\diff^{\nabla,\eta}$, depending on the formulation), one of the pieces must change the associated connection to $\bar\nabla=\nabla+\mathcal{T}^\nabla$, where $\mathcal{T}^\nabla$ is the torsion tensor of $\nabla$, for instance take $\ex_{\bar\nabla}$ together with $\Phi^\nabla$ (or $\diff^{\nabla,\eta}$).

\begin{theorem}
Let $\nabla$ be an arbitrary linear connection on $M$. The exchange map $\ex_\nabla\colon J^1\pi_1\to J^1\pi_1$ given by the local expression \eqref{eq:exchange.map.coord} is well defined. Moreover, we have that:
\begin{enumerate}
\item The exchange map $\ex_\nabla$ is an isomorphism of affine bundles over the identity of $J^1\pi$ that exchanges both affine structures of $J^1\pi_1$ (Figures \ref{fig:iterated.jet.bundle} and \ref{fig:exchange.map}) and whose inverse is $\ex_\nabla^{-1} = \ex_{\nabla+\mathcal{T}^\nabla}$, where $\mathcal{T}^\nabla$ is the torsion tensor of $\nabla$.
\item The exchange map $\ex_\nabla$ is involutive, \ie\ $\ex_\nabla^2=\id_{J^1\pi_1}$, if and only if the connection $\nabla$ is symmetric.
\end{enumerate}
\end{theorem}

\begin{proof}
As in \cite{Mdgn89}, we proof the main assertion by a change of coordinates. Let
\[ (x^i,u^\alpha,u^\alpha_i,\bar u^\alpha_i,u^\alpha_{ii'}) \quand (y^j,v^\beta,v^\beta_j,\bar v^\beta_j,v^\beta_{jj'}) \]
denote a couple of adapted coordinate systems on $J^1\pi_1$. If $\chr^{i''}_{ii'}$ and $\chr^{j''}_{jj'}$ are the Christoffel symbols of $\nabla$ with respect to $(x^i)$ and $(y^j)$ respectively, then
\[ \chr^{j''}_{jj'} = \rbra \frac{\partial^2x^{i''}}{\partial y^j\partial y^{j'}} + \pp[x^i]{y^j}\pp[x^{i'}]{y^{j'}}\chr^{i''}_{ii'} \rket \pp[y^{j''}]{x^{i''}} \,. \]
Besides, we have
\begin{align*}
v^\beta_j =& \rbra \pp[v^\beta]{x^i} + u^\alpha_i\pp[v^\beta]{u^\alpha} \rket \pp[x^i]{y^j} \,,
\intertext{therefore}
\bar v^\beta_j =& \rbra \pp[v^\beta]{x^i} + \bar u^\alpha_i\pp[v^\beta]{u^\alpha} \rket \pp[x^i]{y^j}
\intertext{and}
v^\beta_{jj'} =& \rbra \pp[v^\beta_j]{x^{i'}} + \bar u^\alpha_{i'}\pp[v^\beta_j]{u^\alpha} + u^\alpha_{ii'}\pp[v^\beta_j]{u^\alpha_i} \rket \pp[x^{i'}]{y^{j'}}\\
=& \rbra \frac{\partial^2v^\beta}{\partial x^{i'}\partial x^i}
   + u^\alpha_i\frac{\partial^2v^\beta}{\partial x^{i'}\partial u^\alpha}
   + \bar u^{\alpha'}_{i'}\frac{\partial^2v^\beta}{\partial u^{\alpha'}\partial x^i}
   + u^\alpha_i\bar u^{\alpha'}_{i'}\frac{\partial^2v^\beta}{\partial u^{\alpha'}\partial u^\alpha}
   \nket\\
 & \nbra\
   +\, u^\alpha_{ii'}\pp[v^\beta]{u^\alpha} \rket \pp[x^i]{y^j}\pp[x^{i'}]{y^{j'}}
   + \rbra \pp[v^\beta]{x^i} + u^\alpha_i\pp[v^\beta]{u^\alpha} \rket \frac{\partial^2x^i}{\partial y^{j'}\partial y^j} \,.
\end{align*}
Note also that
\[ \pp[y^j]{x^i}\pp[x^i]{y^{j'}} = \delta^j_{j'} \,,\quad
   \pp[v^\beta]{x^i}\pp[x^i]{y^j} + \pp[v^\beta]{u^\alpha}\pp[u^\alpha]{y^j} = 0  \quand
   \pp[v^\beta]{u^\alpha}\pp[u^\alpha]{v^{\beta'}} = \delta^\beta_{\beta'} \,.\]

Now, let 
\[ (x^i,u^\alpha,\mu^\alpha_i,\bar \mu^\alpha_i,\mu^\alpha_{ii'}) \quand
   (y^j,v^\beta ,\nu^\beta_j ,\bar \nu^\beta_j ,\nu^\beta_{jj'})\]
denote the previous systems of coordinates too, but seen from the image of $\ex_\nabla$. If we assume that
\[ (x^i,u^\alpha,\mu^\alpha_i,\bar \mu^\alpha_i,\mu^\alpha_{ii'}) = 
   (x^i,u^\alpha,  u^\alpha_i,\bar   u^\alpha_i,  u^\alpha_{ii'})\circ\ex_\nabla =
   (x^i,u^\alpha,  \bar u^\alpha_i,  u^\alpha_i,  u^\alpha_{i'i} + (\bar u^\alpha_{i''}-u^\alpha_{i''})\chr^{i''}_{i'i}) \,,\]
we then easily obtain from the coordinate change formulas that
\[ (\nu^\beta_j ,\bar \nu^\beta_j) =
   (  v^\beta_j ,\bar   v^\beta_j)\circ\ex_\nabla = (\bar v^\beta_j, v^\beta_j) \,. \]
On the other hand, since $\mu^\alpha_i=\bar u^\alpha_i$ and $\bar\mu^\alpha_i=u^\alpha_i$, with a proper index rearrangement, we obtain that
\begin{align*}
\nu^\beta_{jj'} - v^\beta_{j'j}
=& \sbra (\mu^\alpha_i-\bar u^\alpha_i)\frac{\partial^2v^\beta}{\partial x^{i'}\partial u^\alpha}
   + (\bar \mu^{\alpha'}_{i'}-u^{\alpha'}_{i'})\frac{\partial^2v^\beta}{\partial u^{\alpha'}\partial x^i}
   + (\mu^\alpha_i\bar \mu^{\alpha'}_{i'}-u^{\alpha'}_{i'}\bar u^\alpha_i)\frac{\partial^2v^\beta}{\partial u^{\alpha'}\partial u^\alpha}
  \nket\\
 &\nbra\
   +\, (\mu^\alpha_{ii'}-u^\alpha_{i'i})\pp[v^\beta]{u^\alpha} \sket \pp[x^i]{y^j}\pp[x^{i'}]{y^{j'}}
   + (\mu^\alpha_i-u^\alpha_i)\pp[v^\beta]{u^\alpha}\frac{\partial^2x^i}{\partial y^{j'}\partial y^j}\\
=& (\mu^\alpha_{ii'}-u^\alpha_{i'i})\pp[v^\beta]{u^\alpha} \pp[x^i]{y^j}\pp[x^{i'}]{y^{j'}}
   + (\bar u^\alpha_{i''}-u^\alpha_{i''})\pp[v^\beta]{u^\alpha}\frac{\partial^2x^{i''}}{\partial y^{j'}\partial y^j}
\end{align*}
and
\[ (\bar v^\beta_{j''}-v^\beta_{j''})\chr^{j''}_{jj'} = (\bar u^\alpha_{i''}-u^\alpha_{i''})\pp[v^\beta]{u^\alpha}
   \rbra \frac{\partial^2x^{i''}}{\partial y^j\partial y^{j'}} + \pp[x^i]{y^j}\pp[x^{i'}]{y^{j'}}\chr^{i''}_{ii'} \rket \,. \]
Subtracting the last two equations and using our initial assumption, we finally obtain that
\[ \nu^\beta_{jj'} = v^\beta_{j'j} + (\bar v^\beta_{j''}-v^\beta_{j''})\chr^{j''}_{jj'} = v^\beta_{jj'}\circ\ex_\nabla \,, \]
which ends the proof of the main statement and shows that $\ex_\nabla$ does not depend on the chosen coordinates and it is well defined.

The two remaining statements are obvious from here.
\end{proof}

% BIBLIOGRAPHY -----------------------------------------------------------------
% Bibliographic styles: amsalpha, amsplain, ieeetr, siam
\bibliographystyle{siam}
\bibliography{t3cft}
% When the bibliography is generated, CAUTION with the alphabetization of the
% "von" particles (which should not be considered).
%\input{t3cft.bbl}
\end{document}